\magnification 1200
\hoffset 0.5 true cm
\hsize 15 true cm
\input amssym.def
\font \tit=cmbx6 scaled \magstep4
\font\bu=cmss8
\def\ZZ{\Bbb{Z}}
\def\QQ{\Bbb{Q}}

\def\LL{\Bbb{L}}
\def\CC{\Bbb{C}}
\def\GG{\Bbb{G}}
\def\BB{\Bbb{B}}
\def\PP{\Bbb{P}}
\def\G{{\cal G}_{\rm mot}}
\def\L{{\rm Lie}\,}
\def\M{\widetilde M}
\def\Hom{\rm Hom}
\def\End{\rm End}
\def\Ext{\rm Ext}

\def\Gr{\rm Gr}
\def\rg{\rm rg}
\def\pn{\par\noindent}
\def\du{\hbox{\bu v}}
\def\ok{\overline k}
\def\gal{{\rm Gal}(\ok/k)}

\input diagrams

\null\vskip+2 true cm

\pn {\centerline {\tit Le radical unipotent du groupe }}
\pn {\centerline {\tit  de Galois motivique d'un 1-motif}}
\vskip 1 true cm

\pn 
{\centerline {Cristiana Bertolin}}
\footnote{}{\pn {\it Math. classification:} 14A99 - 18A99 - 11G99}
\vskip 2 true cm

\par Introduction
\par 1. Rappels sur les 1-motifs et leur groupe de Galois motivique
\par 2. Les $\Sigma-$torseurs et les alg\`ebres de Lie 
\par 3. Etude de ${\Gr}_*^W(\L \G (M))$
\par 4. Preuve du th\'eor\`eme principal
\par Bibliographie
\vskip 1 true cm

\pn  {\bf Introduction}
\vskip 0.5 true cm

\par Soient $M$ un 1-motif defini sur un corps $k$ de
 ca\-ract\'e\-ris\-ti\-que nulle et $\langle M \rangle^\otimes$
la cat\'egorie tannakienne engendr\'ee par
 $M$ (dans une cat\'egorie des syst\`emes de r\'ealisations). Le produit 
tensoriel de la cat\'egorie $\langle M \rangle^\otimes$ permet de d\'efinir
la notion d'alg\`ebre de Hopf commutative dans la cat\'egorie 
${\rm Ind}\langle M \rangle^\otimes $ des Ind-objets de 
$\langle M \rangle^\otimes $. La cat\'egorie des sch\'emas en groupes affines motiviques est l'oppos\'ee de la cat\'egorie 
des alg\`ebres de Hopf commutatives de ${\rm Ind}\langle M \rangle^\otimes.$

 \par Le groupe de Galois motivique $\G (M)$ de $M$ est
le groupe fondamental de la cat\'egorie tannakienne $\langle M \rangle^\otimes$
 engendr\'ee par $M$, i.e.
 le sch\'ema en groupes affine motivique ${\rm Sp}( \Lambda)$, o\`u
 $ \Lambda$ est l'\'el\'ement de $\langle M \rangle^\otimes$ 
universel pour 
la propri\'et\'e suivante:
 pour tout $X$ dans $\langle M \rangle^\otimes,$ il existe un morphisme 
$\lambda_X:  X^{\du} \otimes X \longrightarrow   \Lambda $
 fonctoriel en $X$. Ces morphismes $\{ \lambda_X \}$, qui peuvent se re\'ecrire sous la forme $ X \longrightarrow X \otimes \Lambda$,
 d\'efinissent une action du groupe 
$\G (M)$ sur chaque 
\'el\'ement $X$ de $\langle M \rangle^\otimes$. 
Puisque la d\'efinition de groupe de Galois motivique ne d\'epend que de la cat\'egorie tannakienne engendr\'ee par le 1-motif $M$ et puisque 
deux 1-motifs isog\`enes engendrent la m\^eme cat\'egorie tannakienne,
le groupe $\G (M)$ est d\'efinie \`a 
isog\'enie pr\`es.

\par La filtration par le poids 
$W_{*}$ du 1-motif $M$
induit une filtration croissante $W_*$ en trois crans sur $\G (M)$. En particulier, 
le cran $W_{-1}( \G (M))$ est le radical unipotent du groupe $\G (M).$
Le r\'esultat principal de cet article est un th\'eor\`eme de P. Deligne (th\'eor\`eme  principal 0.1 infra) 
qui affirme que {\it le radical unipotent de l'alg\`ebre de Lie de $\G (M)$
est la vari\'et\'e semi-ab\'elienne d\'efinie par l'action adjointe de l'alg\`ebre de Lie $\Gr_{*}^W(W_{-1}(\L \G (M)))$ sur elle-m\^eme.}
P. Deligne a esquiss\'e la preuve de ce th\'eor\`eme 
dans une lettre qu'il a envoy\'e \`a l'auteur \`a la suite 
d'une question de celle-ci, \`a propos
 du radical unipotent du groupe de Mumford-Tate des 1-motifs. En suivant ses indications (cf. [D01]), dans cet article on se propose 
de donner une preuve compl\`ete de ce r\'esultat.

\par L'id\'ee de la d\'emonstration de ce th\'eor\`eme principal est la suivante: Soit 

$$E=W_{-1}({\underline {\End}}
({\Gr}^W_*(M)))$$

\pn le 1-motif scind\'e de poids -1 et -2
 obtenu \`a partir des endomorphismes du gradu\'e ${\Gr}_*^W (M)$
 de $M$. L'antisym\'etris\'e du produit 

$$P: E \otimes E \longrightarrow E$$

\pn d\'efini par la composition des endomorphismes, munit $E$ d'une
structure d'alg\`ebre de Lie, $(E,[\, ,\,]).$
L'action de $E$ sur le 
gradu\'e $ {\Gr}_*^W (M)$ est d\'ecrite par un morphisme 

$$E \otimes {\Gr}_*^W (M) \longrightarrow {\Gr}_*^W (M)$$

\pn qui munit 
le 1-motif $ {\Gr}_*^W (M)$ d'une structure de $(E,[\, ,\,])$-module de Lie.
\pn La donn\'ee d'une structure de Lie sur le 1-motif $E$ 
est \'equivalente \`a la donn\'ee
d'un $\Sigma$-torseur $d^*{\cal B}$
 sur la vari\'et\'e ab\'elienne $\Gr_{-1}^W(E)$.
Gr\^ace \`a ce $\Sigma$-torseur $d^*{\cal B}$, on peut caract\'eriser le 1-motif $M$:
se donner $M$ est \'equivalent
 \`a se donner le gradu\'e 
${\Gr}_*^W (M)$, un point $k$-rationnel $b$ sur la vari\'et\'e ab\'elienne
 $\Gr_{-1}^W(E)$, et un point $k$-rationnel $\widetilde b$
 dans la fibre du 
$\Sigma$-torseur $d^*{\cal B}$ au dessus du point $b$. 
\pn Soient $B$ la sous-vari\'et\'e ab\'elienne de $\Gr_{-1}^W(E)$ et 
$Z(1)$ le sous-tore de $\Gr_{-2}^W(E)$ tels que 
$(B,Z(1),[\, ,\,])$ est
la plus petite sous-alg\`ebre de Lie de $(E,[\, ,\,])$ dans laquelle vit
 $\widetilde b$. Via le morphisme 
$E \otimes {\Gr}_*^W (M) \longrightarrow {\Gr}_*^W (M)$, cette 
sous-alg\`ebre de Lie agit sur le gradu\'e ${\Gr}_*^W (M)$. De plus,
$(B,Z(1),[\, ,\,])$ est
la plus petite sous-alg\`ebre de Lie de $(E,[\, ,\,])$ qui caract\'erise
la cat\'egorie $\langle M \rangle^\otimes$ au moyen de son action sur  
 ${\Gr}_*^W (M)$: plus pr\'ecisement,
 le foncteur ``prendre le gradu\'e'' 

$${\Gr}^W_*:\langle M \rangle^\otimes \longrightarrow
\langle {\Gr}_*^W (M) \rangle^\otimes$$

\pn induit une \'equivalence de la cat\'egorie
$\langle M \rangle^\otimes$ avec 
 la cat\'egorie des objets de 
\pn $\langle  {\Gr}_*^W (M) \rangle^\otimes$
munis de l'action de l'alg\`ebre de Lie
 $(B,Z(1), [\, , \,])$ (th\'eor\`eme 3.8).
D'un autre c\^ot\'e, si on applique le th\'eor\`eme 8.17 [D90] \`a ce foncteur 
${\Gr}_*^W: \langle M \rangle^\otimes \longrightarrow
 \langle {\Gr}^W_*(M)\rangle^\otimes,$ on obtient que la cat\'egorie
$\langle M \rangle^\otimes$ est \'equivalente \`a la 
cat\'egorie des objets de  $\langle  {\Gr}_*^W (M) \rangle^\otimes$
munis d'une action de $\Gr_{-1}^W(\L \G (M))+\Gr_{-2}^W(\L \G (M)). $
 Gr\^ace \`a la propri\'et\'e universelle du groupe de Galois motivique
$\G (M)$, on conclut que
 $\Gr_{-1}^W(\L \G (M))$
est la vari\'et\'e ab\'elienne $B$ et $W_{-2}(\L \G (M))$ est le tore $Z(1)$
 (corollaire 3.10).
Avec ces notations on peut \'enoncer:
\vskip 0.5 true cm

\pn {\bf 0.1. Th\'eor\`eme principal} 

\pn {\it Le radical unipotent $W_{-1}(\L {\G}(M))$ de l'alg\`ebre de Lie du 
groupe de Galois motivique de $M$ est la vari\'et\'e 
semi-ab\'elienne extension de $B$ par $Z(1)$

$$ 0 \longrightarrow Z(1) \longrightarrow W_{-1}(\L {\G}(M))
 \longrightarrow B  \longrightarrow 0 \leqno(0.1.1)$$

\pn d\'efinie
par l'action adjointe de l'alg\`ebre de Lie $(B,Z(1), [\, , \,])$ sur $B+Z(1).$ }
\vskip 0.5 true cm

\par Une cons\'equence imm\'ediate de ce th\'eor\`eme est 
le calcul de la dimension de l'alg\`ebre de Lie du groupe de Galois motivique ${\G}(M)$
d'un 1-motif:
\vskip 0.5 true cm

\pn {\bf 0.2. Corollaire} 
{\it
$$ \dim \L {\G}(M) = \dim B + \dim Z(1) + \dim \L {\G}(\Gr^W_* M),$$

\pn o\`u on a clairement que $\dim \L {\G}(\Gr^W_* M)=\dim \L {\G}( {\Gr}^W_{-1}M)$ si ${\Gr}^W_{-1}M \not= 0$, et que
  $ \dim \L {\G}(\Gr^W_* M) = 1$ si ${\Gr}^W_{-1}M=0$ et 
 ${\Gr}^W_{-2}M \not= 0$. }
\vskip 0.5 true cm

\par Pour d\'efinir le groupe de Galois motivique d'un 1-motif $M$ et pour pouvoir appliquer le th\'eor\`eme 8.17 de [D90],
on a besoin de savoir que $M$ engendre une cat\'egorie tannakienne. 
Actuellement,
 la seule fa\c con d'obtenir ce r\'esultat est d'identifier la cat\'egorie des 1-motifs \`a une sous-cat\'egorie de la cat\'egorie tannakienne des syst\`emes de r\'ealisations. 
La d\'efinition de groupe de Galois motivique et l'application du 
th\'eor\`eme 8.17 de [D90], sont la raison pour laquelle on est oblig\'e
de supposer que le coprs $k$ est de
 ca\-ract\'e\-ris\-ti\-que nulle (rappelons que la cat\'egorie tannakienne des syst\`emes de r\'ealisations n'est d\'efini que pour $k$ de 
ca\-ract\'e\-ris\-ti\-que nulle).

\par Par contre, la d\'efinition de l'alg\`ebre de Lie $(E,[\,,\,])$ ainsi que
 celle de son action sur le gradu\'e 
$ {\Gr}_*^W (M)$ ne requi\`erent pas l'hypoth\`ese que le corps de d\'efinition de $M$ soit de ca\-ract\'e\-ris\-ti\-que nulle: pour la construction de 
l'alg\`ebre de Lie $(E,[\,,\,])$ et de son action 
$E \otimes {\Gr}_*^W (M) \longrightarrow {\Gr}_*^W (M)$, on se place 
dans la cat\'egorie des 1-motifs (i.e. dans un cadre tout 
 \`a fait g\'eometrique) et non dans la cat\'egorie 
tannakienne des syst\`emes de r\'ealisations. En 
particulier, on n'identifie pas les morphismes 
et les objets de la cat\'egorie des 1-motifs aux 
r\'ealisations correspondantes, mais on d\'efinit
 au fur et \`a mesure  les produits tensoriels et les morphismes
 dont on a besoin: on va essentiellement tensoriser 
des motifs par des motifs purs de poids 0, et comme morphismes on va utiliser 
les projections et les biextensions. Remarquons que {\it les biextentions 
 de deux 1-motifs par un tore 
sont l'interpr\'etation g\'eom\'etrique des
morphismes $M_1 \otimes M_2 \longrightarrow Y(1)$ du produit tensoriel
de deux 1-motifs vers un tore.} Le r\'ef\'er\'e m'a signal\'e qu'on peut
 g\'en\'eraliser cette derni\`ere interpr\'etation pour d\'efinir
 les morphismes   
 $M_1 \otimes M_2 \longrightarrow M_3$ du produit tensoriel
de deux 1-motifs vers un troisi\`eme 1-motif.

\par Pour chaque foncteur fibre $\omega$ de $\langle M \rangle^\otimes$ sur un
 $k$-sch\'ema $S$,
le $S$-sch\'ema $\omega( \G (M)):= {\rm Spec}(\omega (\Lambda)) $ 
est le $S$-sch\'ema en groupes affines
$ {\underline {\rm Aut}}^{\otimes}_S(\omega),$ qui re\-pr\'e\-sen\-te 
le foncteur qui \`a tout $S$-sch\'ema
$T$, $u: T \longrightarrow S$, associe le groupe des automorphismes
 de $\otimes$-foncteurs du foncteur fibre $ u^*\omega.$ En particulier, si 
$\omega_H$
est le foncteur fibre ``r\'ealisation de Hodge'',
 $\omega_H( \G (M))$ est le groupe de Mumford-Tate
 du 1-motif $M$: en d'autres termes, {\it le groupe de Galois motivique de $M$ est l'interpr\'etation g\'eom\'etrique du groupe de Mumford-Tate de $M$}. Ce papier 
compl\`ete donc l'\'etude de la structure du groupe de Mumford-Tate de $M$ que l'auteur
avait entreprise dans \S 1 [B02]: le th\'eor\`eme principal 0.1 
est la version 
compl\`ete et motivique du lemme structural 1.4 [B02]. 

\par La notion de groupe fondamental d'une cat\'egorie tannakienne 
(et donc, en particulier, la notion de groupe de Galois motivique) a \'et\'e introduite par P. Deligne dans [D89] et [D90], o\`u le lecteur trouvera toutes les propri\'et\'es de base de ce groupe. Dans [M94], J. S. Milne \'etudie le groupe de Galois motivique des motifs d\'efinis sur un corps fini. De plus, toujours dans cet article, il d\'emontre que le groupe de Galois motivique des motifs de type CM sur un corps de caract\'eristique nulle, est le groupe de Serre.
 Pour le moment, ceci sont les seules r\'esultats
publi\'es sur les groupes de Galois motiviques. 

\par Dans la premi\`ere section on rappelle rapidement des 
notions de base
sur les 1-motifs d\'efinis sur $k$ et sur leur groupe de Galois motivique.
 Apr\`es un survol sur les $\Sigma$-torseurs, dans la  deuxi\`eme section
 on d\'emontre que la
donn\'ee d'un  $\Sigma$-torseur sur une vari\'et\'e ab\'elienne
est \'equivalente \`a la donn\'ee d'une 
structure d'alg\`ebre de Lie
sur un 1-motif. On termine cette section en 
construisant un exemple 
de $\Sigma$-torseur et en explicitant le crochet de Lie qui lui correspond.
 Dans la section 3 on calcule les gradu\'es 
$\Gr_{-1}^W(\L \G (M))$ et $\Gr_{-2}^W(\L \G (M))$ de l'alg\'ebre de
 Lie du groupe de Galois motivique d'un 1-motif $M$ et
 enfin dans la section 4 on d\'emontre le th\'eor\`eme principal.
\vskip 0.5 true cm

\pn {\it Remerciements}: 
 Je suis tr\`es reconnaissante envers 
 Pierre Deligne pour avoir \'eclair\'e d'un jour nouveau les r\'esultats 
de ma th\`ese de doctorat (cf. [D01]). Je tiens aussi \`a le remercier 
 pour les corrections et les commentaires qu'il a apport\'es  \`a ce texte. Je remercie Lawrence Breen pour avoir mis 
en \'evidence le r\^ole des $\Sigma$-torseurs dans la lettre
 de P. Deligne, et Daniel Bertrand, Uwe Jannsen et Jean-Pierre Wintenberger 
 pour les discussions que nous avons eues sur les 1-motifs. 
 Je tiens \`a remercier aussi le r\'ef\'er\'e pour ses commentaires
 et ses questions.
Enfin, je suis tr\`es reconnaissante envers Michel Waldschmidt qui n'a jamais cess\'e 
de m'encourager et de me soutenir dans ma recherche. 
\pn J'ai 
r\'edig\'e ce papier pendant un s\'ejour \`a l'universit\'e de
M\"unster et \`a l'uni\-ver\-si\-t\'e de
Strasbourg: je tiens \`a les remercier pour leur hospitalit\'e.
\vskip 0.5 true cm

\pn Dans cette article, $k$ est un corps de caract\'eristique nulle et de
cl\^oture alg\'ebrique $\ok,$ et $\gal$ est le groupe de Galois 
de $\ok$ sur $k$.

\vskip 0.5 true cm

\pn {\bf 1. Rappels sur les 1-motifs et leur groupe de Galois motivique}
\vskip 0.5 true cm

\pn 1.1. D'apr\`es [D75] (10.1.10), {\it un 1-motif $M=[X {\buildrel u \over \longrightarrow} G]$ sur $k$} consiste en
\vskip 0.2 true cm

\par (a) un sch\'ema en groupes $X$ sur $k$, 
qui pour la topologie \'etale est un 
sch\'ema en groupes constants d\'efini par un 
${\Bbb Z}$-module libre de type fini;
\vskip 0.2 true cm
\par (b) une vari\'et\'e semi-ab\'elienne d\'efinie sur $k$, i.e.
une extension d'une vari\'
et\' e ab\' elienne $A$ par un tore $Y (1)$ de groupe de cocaract\`eres $Y$,
\vskip 0.2 true cm
\par (c) un morphisme $u:X \longrightarrow G$ de sch\'emas en groupes sur $k.$
\vskip 0.2 true cm

\par
 On doit penser \`a $X$ comme au groupe des caract\`eres
 d'un tore d\'efini sur $k$, i.e. comme \`a
un $\gal$-module de type fini. On identifie $X(\ok)$ \`a un
 $\ZZ$-module libre de type fini, qu'on note $\ZZ^{{\rg}X}.$
La donn\'ee (c) est \'equivalente \`a la donn\'ee
de l'homomorphisme 

$$u:  X(\ok) \longrightarrow G(\ok)$$

\pn $\gal$-\'equivariant.
On note $X$ (resp. $A$, $Y(1)$) le 1-motif $[X  \longrightarrow 0]$ (resp.
$[0  \longrightarrow A]$, $[0  \longrightarrow Y(1)]$).

\par Une isog\'enie entre deux 1-motifs $f=(f_X,f_G): [X_{1} {\buildrel u_{1} \over \longrightarrow} G_{1}]\longrightarrow
[X_{2} {\buildrel u_{2} \over \longrightarrow} G_{2}]$
est  un morphisme de
complexes de sch\' emas en groupes, tel que 
 $f_{X}:X_{1} \longrightarrow X_{2}$ est injectif et de conoyau fini, et
 $f_{G}:G_{1} \longrightarrow G_{2}$ est surjectif et de noyau fini.
On d\' efinit sur $M
=[X {\buildrel u \over \longrightarrow} G]$ 
une filtration croissante $W_{*}$, dite {\it filtration par le poids}:

$$\eqalign{ W_{0}(M)&=M,\cr
 W_{-1}(M)&=[0 \longrightarrow G],\cr
    W_{-2}(M)&=Y(1) .\cr}$$

\pn Si on pose ${\rm
Gr}_{n}^{W}\, {\buildrel {\rm d\acute ef } \over =}\,
W_{n} / W_{n-1},$ on a  ${\rm Gr}_{0}^{W}(M)= X
, {\rm Gr}_{-1}^{W}(M)=  A
$ et $ {\rm Gr}_{-2}^{W}(M)= Y(1).$

\par La notion de biextension permet de donner une description plus sym\'etrique des 1-motifs: consid\'erons le 7-uplet $(X,X^{*},A,A^{*}, v ,v^{*},\psi)$ o\`u

\vskip 0.2 true cm
\par - $X$ et $X^{*}$ sont deux sch\'emas en groupes sur $k$, 
qui pour la topologie \'etale sont des  
sch\'emas en groupes constants d\'efinis par des 
${\Bbb Z}$-modules libres de type fini,
\vskip 0.2 true cm

\par - $A$ et $A^{*}$ sont deux vari\'et\'es ab\'eliennes en dualit\'e
d\'efinies sur $k$,
\vskip 0.2 true cm

\par - $v:X \longrightarrow A$ et 
$v^{*}:X^{*} \longrightarrow A^{*}$ sont deux morphismes 
de sch\'emas en groupes sur $k$ (ou, ce qui est \'equivalent,
$v:X(\ok) \longrightarrow A(\ok)$ et 
$v^{*}:X^{*}(\ok) \longrightarrow A^{*}(\ok)$ sont deux homomorphismes $\gal$-\'equivariants), et 
\vskip 0.2 true cm

\par - $\psi$ est une trivialisation $\gal$-\'equivariante de l'image r\'eciproque par $(v,v^{*})$ de la biextension de Poincar\' e de $(A,A^{*})$.
\vskip 0.2 true cm

\pn A partir de ce 7-uplet $(X,X^{*},A,A^{*}, v ,v^{*},\psi)$ on peut reconstruire le 1-motif $M=[X {\buildrel u \over \longrightarrow} G]$:
$X^*$ s'identifie au groupe des caract\`eres du tore $Y (1)$, 
$v^{*}$ d\'efinit l'exten\-sion $G$ de $A$ par $Y (1),$
et $ v$ et $ \psi$ fournissent l'homomorphisme $u$ 
(cf. [R94] 2.4.1).
\vskip 0.5 true cm

 \pn 1.2. Soient ${\cal M}_{\leq 1}(k)$ la cat\'egorie des 1-motifs sur $k$ 
et $SR(k)$ la cat\'egorie tannakienne des syst\`emes de r\'ealisations 
(pour les cycles de Hodge absolus) sur $k$ d\'efinie par Jannsen dans [J90] I 2.1.
D'apr\`es [D75] \S 10 et [R94] 3.1 on a le foncteur pleinement fid\`ele

$$\eqalign{{\cal M}_{\leq 1}(k) & \longrightarrow SR(k)\cr
 M & \longmapsto ({\rm T}_\sigma(M),{\rm T}_{\rm dR}(M),{\rm T}_{\ell}(M),
I_{\sigma, {\rm dR}}, I_{{\overline \sigma}, {\ell}} )_{{\sigma:k \rightarrow \CC, ~ {\overline \sigma}:\ok \rightarrow \CC} \atop
{{\ell}~{\rm nombre ~ premier}}} \cr}$$

\pn qui \`a chaque 1-motif $M$ associe ses r\'ealisations de Hodge
 ${\rm T}_\sigma(M)$, de de Rham ${\rm T}_{\rm dR}(M)$, $\ell$-adique 
${\rm T}_{\ell}(M)$, et les isomorphismes de comparaison
 $I_{\sigma, {\rm dR}}:{\rm T}_\sigma(M) \otimes_\QQ \CC \longrightarrow 
{\rm T}_{\rm dR}(M)\otimes_k \CC $ et $I_{{\overline \sigma}, {\ell}}:{\rm T}_\sigma(M) \otimes_\QQ 
\QQ_\ell \longrightarrow {\rm T}_{\ell}(M)$.
 On peut donc identifier la cat\'egorie ${\cal M}_{\leq 1}(k)$ 
 \`a une sous-cat\'egorie de $SR(k)$, qu'on notera encore
${\cal M}_{\leq 1}(k)$.
Soit ${\cal M}(k)$ la sous-cat\'egorie tannakienne de la cat\'egorie
$SR(k)$ engendr\'ee par ${\cal M}_{\leq 1}(k)$.
L'objet unit\'e de ${\cal M}(k)$ est le 1-motif $\ZZ(0)=[\ZZ  \longrightarrow 0]$.
On note $M^{\du} \cong {\underline {\Hom}}(M,\ZZ(0)) $
le dual du 1-motif $M$ et $ev_M: M \otimes M^{\du} \longrightarrow \ZZ(0)$ et
$\delta_M: \ZZ(0) \longrightarrow M^{\du} \otimes M$ les morphismes de 
 ${\cal M}(k)$ qui le caract\'erisent (cf. [D90] (2.1.2)).
 Le dual de Cartier de $M$ est le 1-motif

$$M^*=M^{\du} \otimes \ZZ(1). \leqno(1.2.1)$$

\pn {\it La sous-cat\'egorie tannakienne $\langle
M \rangle^\otimes$ engendr\'ee par $M$} est la sous-cat\'egorie pleine de ${\cal M}(k)$ 
d'objets les sous-quotients des sommes de $M^{\otimes\, n} \otimes M^{\du \,\otimes\, m}$
et de foncteur fibre
la restriction du foncteur fibre de ${\cal M}(k)$ \`a $\langle M \rangle^\otimes.$
\par D'apr\`es [By83] (2.2.5), pour chaque plongement $\sigma:k \longrightarrow \CC$
 de $k$ dans $\CC$, le foncteur 
correspondant de la cat\'egorie $\langle M \rangle^\otimes$ vers
 la cat\'egorie des structures de Hodge mixtes est pleinement fid\`ele.
\vskip 0.5 true cm

\pn 1.3. Dans la cat\'egorie ${\cal M}(k)$, un {\it morphisme 
$M_1 \otimes M_2 \longrightarrow Y(1)$} du produit tensoriel 
de deux 1-motifs vers un tore
 est une classe d'isomorphismes de biextensions de $(M_1,M_2)$ par $Y(1)$
(cf. [D75] (10.2.1)).
On pose

$${\Hom}_{{\cal M}(k)}(M_1 \otimes M_2, Y(1))= {\rm Biext}^1(M_1, M_2; Y(1))
\leqno(1.3.1)
$$

\pn o\`u ${\rm Biext}^1(M_1, M_2; Y(1))$ est le groupe des classes 
d'isomorphismes de biextensions de $(M_1,M_2)$ par $Y(1)$.
 En particulier, la classe de 
la biextension de Poincar\'e $\cal P$ de $(A,A^*)$ par $\ZZ(1)$
est l'accouplement de Weil 

$$P_{\cal P}: A \otimes A^* \longrightarrow \ZZ(1).$$

\pn Rappelons que dans [D75] \S 10.2, sous l'hypoth\`ese que $k$ est alg\'ebriquement clos
 Deligne prouve que cette d\'efinition est compatible 
aux r\'ealisations de Hodge, $\ell$-adique et de de Rham. 

\par Soit $M$ un 1-motif de ${\cal M}(k)$.
{\it Une biextension sym\'etrique} $(B,\xi_B) $ de $(M,M)$ par $Y(1)$ 
 est une biextension $B$  de $(M,M)$ par $Y(1)$ (cf. [D75] (10.2.1)) munie d'un 
morphisme de biextensions 
$ \xi_B: s^*B \longrightarrow B,$ o\`u $s^*B$ est l'image inverse de $B$ par le morphisme 
 $s:M \times M \longrightarrow M \times M$ qui permute les facteurs,
tel que la restriction  $d^* \xi_B$ de $ \xi_B$
par le morphisme diagonal $d:M \longrightarrow M \times M$ co\"\i ncide avec 
l'isomorphisme 

$$\nu_B : d^*s^* B \longrightarrow d^*B\leqno(1.3.2)$$

\pn d\'ecoulant de l'identit\'e $s\circ  d =d.$ Le morphisme 
$\xi_B$ est involutif, dans le sens que le morphisme compos\'e $\xi_B \circ s^*\xi_B:
s^*s^* B \longrightarrow s^*B \longrightarrow B$ s'identifie
 \`a l'identit\'e de $B$ (cf. [Br83] 1.7).
{\it La sym\'e\-tri\-s\'ee} d'une biextension $B$ de $(M,M)$ par $Y(1)$
est la biextension 
sym\'etrique  $(B \wedge s^* B, \xi_{B \wedge s^* B})$, o\`u le morphisme $ \xi_{B \wedge s^* B}$
 est donn\'ee canoniquement par le morphisme compos\'e

$$\xi_{B \wedge s^* B}: s^* B \wedge s^* s^* B \longrightarrow s^* B \wedge B
{\buildrel \tau \over \longrightarrow} B \wedge s^* B\leqno(1.3.3)$$

\pn o\`u la premi\`ere fl\`eche 
provient de la relation $ s \circ   s = {\rm id}$ alors que la seconde est 
le morphisme 
$\tau:s^* B \wedge B \longrightarrow B \wedge s^* B$ qui permute les facteurs 
du produit contract\'e. 

\par {\it Un morphisme anti\-sy\-m\'e\-trique $M \otimes M \longrightarrow Y(1)$}
 est une classe d'i\-so\-mor\-phis\-mes de biextensions sym\'etriques de $(M,M)$ par $Y(1)$.
\vskip 0.5 true cm

\pn 1.4. Remarque: 
\par (1) Antisym\'etriser un morphisme $M \otimes M \longrightarrow Y(1)$
\'equivaut donc \`a sym\'etriser la biextension correspondante.
\par (2) Soient $B$ une biextension de $(M_1,M_2)$ par $Y(1)$ et 
 $b: M_1 \otimes M_2 \longrightarrow  Y(1)$ le morphisme qui lui correspond
 via (1.3.1). D'apr\`es une g\'en\'eralisation aux 1-motifs de [D75] (10.2.4), 
l'image inverse $s^* B$ de $ B$ par le 
morphisme $s:M_2 \otimes M_1 \longrightarrow M_1 \otimes M_2$ 
qui permute les facteurs, correspond
au morphisme

$$- b \circ s:M_2 \otimes M_1 \longrightarrow Y(1).$$
\vskip 0.5 true cm

\pn 1.5. {\it La cat\'egorie des sch\'emas affines motiviques 
en $\langle M \rangle^\otimes$} est la cat\'egorie 
oppos\'ee de celle des anneaux commutatifs \`a unit\'e de la cat\'egorie
${\rm Ind}\langle M \rangle^\otimes $ des Ind-objets de 
$\langle M \rangle^\otimes $ (cf [D89] \S 5). L'objet final de cette cat\'egorie 
est le sch\'ema affine motivique ${\rm Sp}(\ZZ(0))$  
 d\'efini par l'anneau $\ZZ(0)$ de
 $\langle M \rangle^\otimes $ (rappelons qu'il y a un foncteur canonique pleinement fid\`ele de la cat\'egorie $\langle M \rangle^\otimes $
vers la cat\'egorie ${\rm Ind}\langle M \rangle^\otimes $).
 Un sch\'ema en groupes affine motivique est un objet en groupes de la
 cat\'egorie des sch\'emas affines motiviques. 
L'alg\`ebre de Lie d'un sch\'ema en groupes affine motivique
est un pro-objet $\rm L$ de $\langle M \rangle^\otimes$
muni d'une structure d'alg\`ebre de Lie, i.e. $\rm L$ est muni
d'une application antisym\'etrique
$[\, , \,]: {\rm L} \otimes {\rm L} \longrightarrow {\rm L}$
qui satisfait l'identit\'e de Jacobi.

\par {\it Le groupe de Galois motivique ${\G}(M)$}
 d'un 1-motif $M$ de ${\cal M}(k)$ est
 le sch\'ema en groupes affine motivique ${\rm Sp}(\Lambda),$ o\`u
 $ \Lambda$ est l'\'el\'ement de $\langle M \rangle^\otimes$ universel pour 
la propri\'et\'e suivante: pour tout $X$ dans $\langle M \rangle^\otimes$ il existe 
un morphisme 
$$\lambda_X: X^{\du} \otimes X \longrightarrow \Lambda \leqno(1.5.1)$$
\pn fonctoriel en $X$, i.e.
tel que pour tout $f:X \longrightarrow Y$ dans $\langle M \rangle^\otimes$ le diagramme 

$$\matrix{Y^{\du} \otimes X&{\buildrel f^t \otimes 1 \over \longrightarrow}&
X^{\du} \otimes X\cr
     {\scriptstyle 1 \otimes f}\downarrow~~~~ && ~~~~~ \downarrow
 {\scriptstyle \lambda_X}  \cr
   Y^{\du} \otimes Y& {\buildrel \lambda_Y \over \longrightarrow}&\Lambda \cr}$$

\pn soit commutatif. La propri\'et\'e universelle de $\Lambda$ est que 
pour tout $U$ dans ${\rm Ind}\langle M \rangle^\otimes$, l'application

 $$\eqalign{{\Hom}(\Lambda, U)& \longrightarrow \{ u_X:
  X^{\du} \otimes X \longrightarrow U, ~~ {\rm {fonctoriels~ en~}} X \} \cr
f & \longmapsto  f  \circ \lambda_X \cr }\leqno(1.5.2)$$
 
\pn est une bijection.  
L'existence de ${\G}(M)$ est demontr\'e dans [D90] 8.4, 8.10, 8.11 (iii).
Puisque le groupe ${\G}(M)$ ne d\'epend que de la cat\'egorie tannakienne 
engendr\'e par le 1-motif $M$, la notion de groupe de Galois 
motivique est stable par dualit\'e et par isog\`enies. 

\par Les morphismes (1.5.1), qui peuvent se r\'eecrire sous la forme $ X \longrightarrow  X \otimes \Lambda,$
d\'efinissent une action du groupe ${\G}(M)$ sur chaque 
$X$ de $\langle M \rangle^\otimes$. En particulier, le morphisme 
$ \Lambda \longrightarrow  \Lambda \otimes \Lambda$
repr\'esente l'action de  ${\G}(M)$ sur lui-m\^eme par automorphismes int\'erieurs
(cf. [D89] 6.1).

\par La filtration par le poids de $M$ induit sur 
le groupe de Galois motivique ${\G}(M)$
une filtration croissante $W_*$ d\'efinie sur $k$ (cf. [S72] Chapitre 2 \S 2):

\vskip 0.5 true cm
\pn $ \quad  W_{0}({\G}(M))={\G}(M) $
\vskip 0.3 true cm
\pn $ \quad  W_{-1}({\G}(M))=\{ g \in {\G}(M) \, \, \vert \, \,
(g - id)M \subseteq   W_{-1}(M)  ,(g - id) W_{-1}(M)
\subseteq W_{-2}(M), (g - id) W_{-2}(M)=0 \} , $
\vskip 0.3 true cm
\pn $  \quad  W_{-2}({\G}(M))=\{ g \in {\G}(M) \, \, \vert \, \,
(g - id) M \subseteq W_{-2}(M),  (g - id)  W_{-1}(M) =0\}, $
\vskip 0.3 true cm
\pn $ \quad  W_{-3}({\G}(M))=0.$
\vskip 0.5 true cm

\pn Il est \'evident que $ W_{-1}({\G}(M))$ est unipotent.
D'apr\`es l'analogue motivique de [By83] \S 2.2, 
${\rm Gr}_{0}^{W}({\G}(M))$ agit trivialement sur ${\rm Gr}_{0}^{W}(M)$,
par homoth\'eties sur ${\rm Gr}_{-2}^{W}(M)$ et son
image dans le groupe des automorphismes de
${\rm Gr}_{-1}^{W}(M)$
est le groupe de Galois motivique de la vari\'et\'e ab\'elienne
 ${\rm Gr}_{-1}^{W}(M). $ Par cons\'equent ${\rm Gr}_{0}^{W}({\G}(M))$
est r\'eductif et $W_{-1}({\G}(M))$ est 
le radical unipotent de ${\G}(M)$.
\vskip 0.5 true cm

\pn 1.6. Si $\omega$ est 
un foncteur fibre de $\langle M \rangle^\otimes$ sur un $k$-sch\'ema $S$,
$\omega (\Lambda)$ est l'alg\'ebre de Hopf dont le spectre ${\rm Spec} (\omega (\Lambda))$ est le $S$-sch\'ema affine 
$ {\underline {\rm Aut}}^{\otimes}_S(\omega)$ qui repr\'esente 
le foncteur qui \`a tout $S$-sch\'ema
$T$, $u: T \longrightarrow S$, associe le groupe des automorphismes
 de $\otimes$-foncteurs du foncteur $ X \longmapsto u^*\omega(X)$ 
de $\langle M \rangle^\otimes$ 
dans les faisceaux localement libres de rang fini sur $T$ (cf. [D90] (8.13.1)).

\par D'apr\`es le formalisme de [D89] 5.11,
le groupe de Galois motivique ${\G}(M)$ correspond \`a la donn\'ee,
 pour tout foncteur fibre $\omega$ sur tout $k$-sch\'ema $S$, du 
$S$-sch\'ema en groupes affine
$ {\underline {\rm Aut}}^{\otimes}_S(\omega)$, de formation 
compatible \`a tout changement de base 
$S'/S$, i.e. 
${\G}(M)$ correspond \`a $\{ \omega \longmapsto
 {\underline {\rm Aut}}^{\otimes}_S(\omega) \}$. 
\vskip 0.5 true cm

\pn 1.7. Exemples: 
\pn (1) Le groupe de Galois motivique de l'objet unit\'e $\ZZ(0)$ 
de ${\cal M}(k)$ est le sch\'ema affine motivique ${\rm Sp}(\ZZ(0))$.
Si $\omega_H$ est le foncteur fibre ``r\'ealisation de Hodge'', on a que 
$ \omega_H \big( \G(\ZZ(0))\big):={\rm Spec}\big(\omega_H (\ZZ(0)) \big)=
{\rm Spec}( \QQ)$, qui est le groupe de Mumford-Tate de $\ZZ(0)$.
\pn (2) L'alg\`ebre de Lie du groupe de Galois motivique du tore $\ZZ(1)$ 
est le tore lui-m\^eme, i.e. $\L \big(\G(\ZZ(1))\big) =\ZZ(1)$.
Si on applique le foncteur fibre ``r\'ealisation de Hodge'', on obtient 
que l'alg\`ebre de Lie du groupe de Mumford-Tate de $\ZZ(1)$ est 
la structure de Hodge de Tate.
\vskip 0.5 true cm

\pn {\bf 2. Les $\Sigma$-torseurs et les alg\`ebres de Lie}
\vskip 0.5 true cm

\pn 2.1. Soient $P$ et $G$ deux groupes ab\'eliens dans un topos quelconque et
 $L$ un $G$-torseur sur $P$. Notons 
 $m: P \times P \longrightarrow P$ la loi de groupe de $P$ et 
 $s: P \times P \longrightarrow P \times P$ le morphisme qui permute les facteurs.
On pose 

$$\Lambda(L)= m^*L \wedge p_1^* L^{-1} \wedge p_2^* L^{-1}.$$ 

\pn  o\`u $p_i: P \times P \longrightarrow P$ d\'esigne la projection canonique
sur le $i$-\`eme facteur. 
C'est un torseur sur $P  \times P$ muni
d'un isomorphisme canonique $\xi_{L}: s^* \Lambda(L) \longrightarrow \Lambda(L)$
qui est d\'efini par la commutativit\'e de la loi de groupe de $P$, et
dont l'image inverse par le morphisme diagonal 
est l'isomorphisme  $ d^*s^* \Lambda(L) \longrightarrow d^*\Lambda(L)$
d\'efini par la relation $s\circ  d =d.$
On note 

$$\theta(L)=
m^*_{123} L  \wedge m^*_{12} L^{-1} \wedge m^*_{13} L^{-1} 
\wedge m^*_{23} L^{-1} \wedge 
p_{1}^* L \wedge p_{2}^* L \wedge p_{3}^* L$$

\pn o\`u $m_{123}: P^3  \longrightarrow P$
 (resp. $ m_{ij}: P^3  \longrightarrow P$) 
est la fl\`eche d\'efinie par $m_{123}(x_1,x_2,x_3)$
\pn $=x_1+x_2+x_3$
(resp. $ m_{ij}(x_1,x_2,x_3)=x_i+x_j$) et $p_i:P \times P \times P
 \longrightarrow P$ d\'esigne la projection canonique sur le 
$i$-\`eme facteur. 
Si on confronte les fibres de $\theta(L)$ et 
$(m \times 1)^* \Lambda(L) \wedge p_{13}^* \Lambda(L)^{-1} \wedge p_{23}^* \Lambda(L)^{-1}$ 
au dessus du point g\'en\'eral $(x,y,z)$ de $P \times P \times P$,
on observe que 
les morphismes contraction ${\underline 0} \longrightarrow L_{z}^{-1}L_{z}$ 
et ${\underline 0} \longrightarrow L_{x}^{-1}L_{x}$
d\'efinissent
les isomorphismes canoniques

$$\eqalign{
\chi_1:& \theta(L) \longrightarrow
(m \times 1)^* \Lambda(L) \wedge p_{13}^* \Lambda(L)^{-1} \wedge p_{23}^* \Lambda(L)^{-1},
\cr
\chi_2:& \theta(L) \longrightarrow
 (1 \times m)^* \Lambda(L) \wedge p_{12}^* \Lambda(L)^{-1}  \wedge p_{13}^*
 \Lambda(L)^{-1}. \cr}$$

\par {\it Un $G$-torseur cubiste $(L,(B, \xi_B),\alpha)$} sur $P$ 
est un $G$-torseur $L$ sur $P$ muni d'une biextension 
sym\'etrique $(B, \xi_B)$ de $(P,P)$ par $G$ et d'un morphisme de torseurs
$\alpha: B \longrightarrow \Lambda(L)$, tel que 
$\xi_{L} \circ  s^*\alpha =\alpha  \circ \xi_B$ et tel que 
le diagramme suivant soit commutatif

$$\matrix{ 
(m \times 1)^* B \wedge p_{13}^* B^{-1} \wedge p_{23}^* B^{-1} & \longrightarrow& 
(m \times 1)^* \Lambda(L) \wedge p_{13}^* \Lambda(L)^{-1} \wedge p_{23}^* \Lambda(L)^{-1}   \cr
 {\scriptstyle \psi_B }    \downarrow &      &  \downarrow {\scriptstyle \psi_{L} } \cr
 (1 \times m)^* B \wedge p_{12}^* B^{-1} \wedge p_{13}^* B^{-1}   & \longrightarrow & 
(1 \times m)^* \Lambda(L)  \wedge p_{12}^* \Lambda(L)^{-1} \wedge p_{13}^* \Lambda(L)^{-1} \cr}$$

\pn o\`u les fl\`eches horizontales sont induites par $\alpha, \psi_{L}=\chi_1^{-1} \circ \chi_2$ et 
$\psi_B = s_2 \circ  s_1^{-1}$ o\`u $s_2$ et $ s_1$ sont les
 sections
qui d\'efinissent les deux lois de composition partielles sur la biextension $B$.
Le morphisme $\alpha:  B \longrightarrow \Lambda(L)$ impose, par 
transfert de structure, 
une structure de biextension \`a $\Lambda(L)$
qui est automatiquement une structure de biextension 
sym\'etrique, \`a cause de l'existence de l'isomorphisme canonique $\xi_{L}: s^*\Lambda(L) \longrightarrow \Lambda(L) $.
\vskip 0.5 true cm

\pn 2.2. 
{\it Un $G$-torseur sym\'etrique} $(L,\lambda)$ sur $P$
est un $G$-torseur $L$ sur $P$ muni d'un
isomorphisme de torseurs $ \lambda: {inv}^* L \longrightarrow L$, o\`u
${inv}:P \longrightarrow P$ est la loi d'inverse de $P$.

\par Puisque ${\rm Biext}^1(P_1,P_2;G)$ est additif en les variables $P_1$
et $P_2$,
si $B$ est une biextension de $(P,P)$ par $G$, on a l'isomorphisme de biextension 

$$
(p_1 \times p_3)^* B \wedge (p_1 \times p_4)^* B
 \wedge (p_2 \times p_3)^* B  \wedge (p_2 \times p_4)^* B
\longrightarrow (m \times m)^* B$$

\pn o\`u $p_i: 
P \times P \times P \times P  \longrightarrow P$ est la projection canonique sur le 
$i$-\`eme facteur. Cet isomorphisme fournit le morphisme
$(p_1 \times p_4)^* B \wedge (p_2 \times p_3)^* B \longrightarrow
(m \times m)^* B \wedge (p_1 \times p_3)^* B^{-1} \wedge
(p_2 \times p_4)^* B^{-1},$
dont l'image inverse par le morphisme diagonal $
P^2  \longrightarrow P^2 \times P^2 $ est le morphisme 

$$\gamma_B: B \wedge s^*B \longrightarrow \Lambda ( d^*B).\leqno(2.2.1)$$

\par {\it Un $\Sigma-G$-torseur $((L,\lambda),(B, \xi_B),\alpha,\beta)$} sur $P$ 
consiste en un $G$-torseur cubiste $(L, (B, \xi_B),\alpha)$ dont le 
$G$-torseur sous-jacent est sym\'etrique $(L, \lambda)$, et en un 
 morphisme de torseurs  
$\beta:d^* B \longrightarrow L^2$ tel que 
le diagramme suivant soit commutatif

$$\matrix{ 
& &B^2 & &  \cr
&{}^{c_B}\swarrow & & \searrow^{\alpha^2} &  \cr
\Lambda(d^* B) & & {\buildrel \Lambda(\beta) \over \longrightarrow}
 & & \Lambda(L^2) \cr}$$

\pn o\`u $c_B: B^2  {\buildrel {1 \wedge \xi_B^{-1}} \over \longrightarrow}
B \wedge s^* B {\buildrel \gamma_B \over \longrightarrow} \Lambda(d^* B)$.
 Le morphisme  $\beta$ exprime 
la compatibilit\'e entre la structure sym\'etrique et la structure cubiste de $L$.
Heuristiquement, ceci s'explique ainsi: 
si $\Delta_1: P^2  \longrightarrow  P^3$ est le morphisme d\'efini par $
\Delta_1(x,y)=(x+y,-x,-y)$, alors il existe un isomorphisme 

$$\Delta_1^* \theta (L) \longrightarrow \Lambda(L) \wedge \Lambda(inv^*L)^{-1}.
\leqno(2.2.2)$$

\pn La structure cubiste sur $L$ munit $\Lambda(L)$ d'une structure de biextension, 
qui par les  applications $\chi_i$ et $\Delta_1$
nous fournit une section de $\Delta_1^* \theta (L)$. L'isomorphisme (2.2.2)
 fait correspondre \`a cette section le morphisme de torseurs $\ell_L: \Lambda(inv^*L)
\longrightarrow \Lambda(L)$.  Dire que la structure cubiste de $L$ 
 est compatible \`a la structure sym\'etrique de $L$ signifie que  $
\ell_L$ co\"\i ncide avec le morphisme $
\Lambda(\lambda): \Lambda(inv^*L) \longrightarrow \Lambda(L)$.
\vskip 0.5 true cm

\pn 2.3. Dans SGA7 I expos\'e VII \S 3, Grothendieck
 prouve que ${\rm Biext}^1(P,P;G) \cong 
{\Ext}^1(P \otimes^{\LL} P,G).$ 
De fa\c con analogue, dans [Br83] \S 8 Breen fournit l'interpr\'etation
homotopique des biextensions sym\'etriques,
des torseurs cubistes et des $\Sigma$-torseurs.

 Gr\^ace \`a cette intepretation, on a l'analogie suivant:
\par - une biextension $B$ de $(P,P)$ par $G$ est l'analogue d'une application $f:P \otimes P \longrightarrow G$;
\par - la biextension sym\'etris\'ee $B \wedge s^*B$ de $B$ est l'analogue de  
l'application sym\'etris\'ee $f(p,p')+f(p',p)$ de $f(p,p')$;
\par - imposer une structure de biextension sur 
$\Lambda(L)= m^*L \wedge p_1^* L^{-1}
\wedge p_2^* L^{-1}$ est analogue \`a imposer \`a la premi\`ere diff\'erence de 
$q:P \longrightarrow G,$
$Q(p,p')= q(p+p')-q(p)-q(p')$ 
d'\^etre bilin\'eaire, i.e. les torseurs cubistes
sont les analogues des applications de degr\'e 2;

\par - un $G$-torseur sym\'etrique $L$ sur $P$ est l'analogue d'une application 
 $q:P \longrightarrow G$ tel que $q(p)=q(-p)$ pour chaque $p$ dans $P$;
\par - imposer 
la compatibilit\'e entre la structure cubiste de $L$ et la structure de sym\'etrie de $L$ est analogue \`a imposer que 
 $q:P \longrightarrow G$ soit de degr\'es 2 et que 
$q(p)=q(-p)$ pour chaque $p$ dans $P,$
i.e. les $\Sigma$-torseurs sont les analogues des applications quadratiques.
\vskip 0.5 true cm

\pn 2.4. Apr\`es ce rappel sur les $\Sigma$-torseurs, on va maintenant
d\'emontrer que la
donn\'ee d'un  $\Sigma$-torseur sur une vari\'et\'e ab\'elienne
est \'equivalente \`a la donn\'ee d'une 
structure d'alg\`ebre de Lie sur un 1-motif.
\par Soit $M=[0 \longrightarrow  G]$ un 1-motif 
de ${\cal M}(k)$, o\`u $W_{-1}(M)=G$ est une extension d'une vari\'et\'e ab\'elienne
$A$ par un tore $Y(1)$. Notons $d:A \longrightarrow A \times A$
le morphisme diagonal de $A$.
\vskip 0.5 true cm

\pn 2.5. {\bf Proposition}
\pn {\it 
La donn\'ee d'un  $\Sigma-Y(1)$-torseur $((L,\lambda),(B, \xi_B),\alpha,\beta)$
sur $A$
\'equivaut \`a la donn\'ee d'une structure d'alg\`ebre de Lie $[\, ,\,]:
M \otimes M \longrightarrow M$ sur $M$. 
\pn Plus pr\'ecisement,
le crochet de Lie $[\, ,\,]$
et la biextension sym\'etrique $(B, \xi_B)$ se d\'efinissent mutuellement via (1.3.1),
 et le $\Sigma-Y(1)$-torseur correspondant au crochet de Lie $[\, ,\,]$ est  
la restriction 
$d^* B$ de la biextension sym\'etrique $B$
. }
\vskip 0.5 true cm

\pn PREUVE:  Soit $((L,\lambda),(B, \xi_B),\alpha,\beta)$ un  $\Sigma-Y(1)$-torseur sur $A$.
Par d\'efinition, la biextension sym\'etrique
$ (B, \xi_{B})$ d\'efinit une application bilin\'eaire
 antisym\'etrique 
$[\, , \,]: A \otimes A \longrightarrow Y(1)$. Cette application se rel\`eve 
en un morphisme $[\, , \,]: G \otimes G \longrightarrow Y(1)$ bilin\'eaire
 antisym\'etrique, qui \`a cause des poids
est la seule composante non nulle d'un morphisme $[\, , \,]: M \otimes M \longrightarrow M$.
 Pour chaque $x$ et $y$ dans $M=W_{-1}(M)$,
 le produit $[x,y]$ est dans $W_{-2}(M)=Y(1)$ et donc le produit $[[x,y],z]$ est nul pour 
chaque $z$ dans $M$. L'application $[\, ,\,]$ v\'erifie
donc trivialement 
l'identit\'e de Jacobi. 
\par Soit $[\, , \,]: M \otimes M \longrightarrow M$ une 
structure d'alg\`ebre de Lie sur $M$. A cause des poids, elle a pour seule composante 
non nulle $[\, , \,]: G \otimes G \longrightarrow Y(1)$. Puisque $Y(1) \otimes Y(1)$ 
est un motif de poids -4, ce dernier morphisme se factorise par
un morphisme 
$ A \otimes A \longrightarrow Y(1)$ toujours bilin\'eaire et antisym\'etrique, et
 qui via (1.3.1) d\'efinit une biextension sym\'etrique
$ (B,\xi_{ B})$ de $(A,A)$ par $Y(1)$. D'apr\`es [Br83] 5.7 
la restriction 
$d^* B$ de $B$ par le morphisme diagonal 
$d:A \longrightarrow A \times A$ est muni d'une 
$\Sigma$-structure canonique  $((d^*B ,(-1)), ( B \wedge s^*B ,
 \xi_{B \wedge s^*B}), \gamma_{ B}, 1 \wedge \nu_{ B}).$
\vskip 0.5 true cm

\pn 2.6. Appliquons cette proposition aux $\Sigma$-torseurs obtenus
\`a partir de biextensions: Soit ${\cal B}$ une biextension de $(A,A)$
par $Y(1)$. La restriction $d^* {\cal B}$ de ${\cal B}$
par le morphisme diagonal $d:A \longrightarrow A \times A$ est 
un $Y(1)$-torseur sym\'etrique 
dont la structure de sym\'etrie est donn\'ee par l'iso\-mor\-phis\-me de torseurs $-1 : (-1)^* d^* {\cal B}
\longrightarrow d^* {\cal B}$, o\`u $-1: A \longrightarrow A$ est la loi 
d'inverse de la vari\'et\'e ab\'elienne 
$A$.  
De plus,
d'apr\`es [Br83] 5.7 ce torseur 
est de fa\c con canonique un $\Sigma-Y(1)$-torseur sur $A:$ 
sa $\Sigma$-structure canonique est donn\'ee par le quadruplet

$$((d^* {\cal B}, -1), ({\cal B} \wedge s^*{\cal B}, \xi_{{\cal B} \wedge s^*{\cal B}}), \gamma_{\cal B}, 1 \wedge \nu_{\cal B})$$
 
\pn o\`u $\gamma_{\cal B}: {\cal B} \wedge s^*{\cal B} 
\longrightarrow \Lambda ( d^*{\cal B})$ et 
$ 1 \wedge \nu_{\cal B} :d^*({\cal B} \wedge s^*{\cal B}) \longrightarrow (d^*{\cal B} )^2$ 
(cf. (1.3.2), (1.3.3) et (2.2.1)).
\vskip 0.5 true cm

\pn 2.7. {\bf Corollaire}
\pn {\it
La donn\'ee du $\Sigma-Y(1)$-torseur $d^* {\cal B}
=((d^* {\cal B},-1), ({\cal B} \wedge s^*{\cal B}, \xi_{{\cal B} \wedge s^*{\cal B}}), \gamma_{\cal B}, 1 \wedge \nu_{\cal B})$   sur $A$
\'equivaut \`a la donn\'ee du crochet de Lie $[\, , \,]:M \otimes M \longrightarrow M$
 correspondant \`a la biextension sym\'etrique 
$({\cal B} \wedge s^*{\cal B}, \xi_{{\cal B} \wedge s^*{\cal B}})$ via 
(1.3.1). }
\vskip 0.5 true cm

\pn PREUVE: Puisque $ s \circ d = d,$ on observe que 
  $d^*( {\cal B} \wedge s^*{\cal B}) = 2\, d^* {\cal B}$.
A partir de cette remarque, il suffit d'appliquer la proposition 2.5 et de se rappeler 
qu'en travaillant modulo isog\`enie, on peut inverser le 2. 
\vskip 0.5 true cm

\pn 2.8. Exemple: Soient $x$ et $y$ deux entiers non nuls.
A partir de la biextension  de
 Poincar\'e ${\cal P}^{*}$ de $(A^*,A)$, on va
construire de fa\c con explicite une biextension
$\PP$ de $(A^x + (A^*)^y,A^x + (A^*)^y)$ par
$\ZZ^{x \cdot y}(1)$, et le crochet de Lie
sur le 1-motif scind\'e $(A^x + (A^*)^y)+\ZZ^{x \cdot y}(1)$
d\'efini par
 le $\Sigma-\ZZ^{x \cdot y}(1)$-torseur 
$d^*\PP$, o\`u $d$ est le morphisme diagonal de 
$A^x + (A^*)^y$.
\par Consid\'erons $( x \cdot y)$-copies de 
la biextension de Poincar\'e ${\cal P}$ de $(A,A^*)$, 
$( x \cdot y)$-copies de la biextension de
 Poincar\'e ${\cal P}^{*}$ de $(A^*,A)$, $ x^2$-copies
 de la biextension trivial ${\underline 0}_{A}= A \times A \times {\GG}_m$
 de $(A,A)$ par ${\GG}_m$ et $y^2$-copies de la biextension 
trivial ${\underline 0}_{A^*}= A^* \times A^* \times {\GG}_m$ de $(A^*,A^*)$
par ${\GG}_m$. Parce que la cat\'egorie
${\rm Biext}^1(A,A^*;{\GG}_m)$ est additive en les variables $A$ et $A^*$
 (voir SGA 7 I Expos\'e VII (2.4.2)), \`a partir de ces $ (  x+y)^2$-biextensions on obtient une biextension $\BB$ de $(A^x + (A^*)^y,A^x + (A^*)^y)$ par
${\GG}_m$ qui a la propri\'et\'e suivante:

$$\eqalign{ 
&{\BB}_{| (0\times \dots 0 \times A \times 0 \dots \times 0, 0 \times \dots 0 \times A^* \times 0 \dots \times 0)}= {\cal P} \cr
&{\BB}_{| (0 \times \dots 0 \times A^* \times 0 \dots \times 0, 0 \times  \dots 0 \times A \times 0 \dots \times 0) }={\cal P}^*\cr
&{\BB}_{|(0 \times \dots 0 \times A^* \times 0 \dots \times 0, 0 \times \dots 0 \times A^* \times 0 \dots \times 0)  }={\underline 0}_{A^*}\cr
&{\BB}_{| (0 \times \dots 0 \times A \times 0 \dots \times 0, 0 \times  \dots 0 \times A \times 0 \dots \times 0) }={\underline 0}_{A}\cr}\leqno(2.8.1)$$

\pn De fa\c con explicite, ${\BB}$ se construit ainsi (cf. SGA7 I expos\'e VII 2.4 ou 
[Br83] (1.2.3)):
soient $pr_i^A: A^x + (A^*)^y \longrightarrow A$ et
 $pr_j^{A^*}: A^x + (A^*)^y \longrightarrow A^*$ les projections canoniques
respectivement  
sur le $i$-\`eme et le $x+j$-\`eme facteur pour $i=1, \dots,x$ et $j=1, \dots,y$.
Les biextensions

$$\eqalign{ 
&{\BB}_1=
\Big( \bigwedge_{i=1}^x (pr_i^A \times id_A)^* {\underline 0}_{A} \Big) \bigwedge 
\Big( \bigwedge_{j=1}^y (pr_j^{A^*} \times id_A)^* {\cal P}^* \Big)\cr
&{\BB}_2=
\Big( \bigwedge_{i=1}^x  (pr_i^A \times id_{A^*})^* {\cal P} \Big) \bigwedge 
\Big( \bigwedge_{j=1}^y  (pr_j^{A^*} \times id_{A^*})^* {\underline 0}_{A^*} \Big)\cr}
\leqno (2.8.2)$$

\pn sont respectivement des biextensions de $(A^x + (A^*)^y,A)$ par
${\GG}_m$ et $(A^x + (A^*)^y,A^*)$ par ${\GG}_m$. Avec ces notations on a 

$${\BB}=
\Big( \bigwedge_{i=1}^x (id_{A^x + (A^*)^y} \times pr_i^A)^*{\BB}_1 \Big) \bigwedge 
\Big( \bigwedge_{j=1}^y (id_{A^x + (A^*)^y} \times pr_j^{A^*} )^*{\BB}_2 \Big).
\leqno(2.8.3)$$

\pn A partir de 
$( x \cdot y)$-copies de la biextension $\BB$,
 en utilisant la propri\'et\'e d'additivit\'e de la cat\'egorie 
${\rm Biext}^1(A^x + (A^*)^y,A^x + (A^*)^y;{\GG}_m)$ cette fois par 
rapport \`a l'argument ${\GG}_m,$ on construit une biextension $ \PP$ de $(A^x + (A^*)^y,A^x + (A^*)^y)$ par 
${\ZZ}^{x \cdot y} (1)$: pour 
$l=1, \dots, x \cdot y$, soient $s_l: {\GG}_m \longrightarrow 
{\ZZ}^{x \cdot y}(1)$ les sections qui d\'ecrivent 
${\ZZ}^{x \cdot y}(1)$ comme somme directe de ses facteurs.
 La biextension $ \PP$ est le produit contract\'e 
des images directes de la biextension $\BB$ par les sections $s_l$:
$${\PP}=
\bigwedge_{l=1}^{x \cdot y} (s_l)_* \, \BB.$$

\par Posons $C=A^x + (A^*)^y, W(1)={\ZZ}^{x \cdot y} (1)$ et notons

$$N=C+W(1)$$

\pn le 1-motif scind\'e de seules composantes non nulles $W_{-1}(N)=C$ et 
$W_{-2}(N)=W(1).$
Soit $d:C \longrightarrow C \times C$
 le morphisme diagonal de $C$. La restriction 
$d^* {\PP}$ de la biextension  ${\PP}$
de $(C,C)$ par $W(1)$, est munie d'une structure de 
 $\Sigma-W(1)$-torseur sur $C$.
 D'apr\`es le corollaire 2.7 
la donn\'ee de ce $\Sigma-W(1)$-torseur est 
\'equivalente \`a la donn\'ee du crochet de Lie $[\, , \,]:N \otimes N \longrightarrow N$
 correspondant \`a la biextension sym\'etrique 
$({\PP} \wedge s^*{\PP}, \xi_{{\PP} \wedge s^*{\PP}})$ via 
(1.3.1). Par construction 
ce crochet de Lie $[\, , \,]:N \otimes N \longrightarrow N$
s'explicite ainsi

$$\eqalign{ 
&[\, , \,]_{|  A^x \otimes A^x}=0\cr
&[\, , \,]_{| (A^*)^y  \otimes (A^*)^y }=0\cr
&[\, , \,]_{|  A^x \otimes (A^*)^y }: A^x 
\otimes (A^{*})^y \longrightarrow {\ZZ}^{x \cdot y}(1) \cr
& (a_1, \dots, a_x), (b_1, \dots, b_y) \longmapsto (<a_i,b_j>)_{i,j}\cr}\leqno(2.8.4)$$

\pn o\`u $<\, ,\,>: A \otimes A^* \longrightarrow {\ZZ}(1)$ est l'accouplement de Weil d\'efini
via (1.3.1) par la biextension de Poincar\'e ${\cal P}$ de $(A,A^*)$, qui est une biextension 
sym\'etrique.
Par cons\'equent le produit $[\, , \,]:N \otimes N \longrightarrow N,$ 
de seule composante non nulle  $[\, , \,]:C \otimes C \longrightarrow W(1),$ s'id\'entifie  
\`a $(x \cdot y)$-copies de 
l'ac\-cou\-ple\-ment de Weil $<\, ,\,>:A \otimes A^* \longrightarrow {\ZZ}(1)$. \vskip 0.5 true cm

\pn 2.9. Remarque: Si on identifie $A$ et $A^*$ avec ${\underline {\Ext}}^1(A^*,{\GG}_m)$ et
 ${\underline {\Ext}}^1(A,{\GG}_m)$ res\-pec\-ti\-ve\-ment,
 un point dans la fibre $  (d^* {\PP})_c $  de $d^*{\PP}$ au dessus d'un point
 $c=(c_1,c_2)$ de $A^x \times (A^*)^y(k)$, correspond 
\`a un 1-motif $N_c$ tel que

\par (i) $W_{-2}(N_c)= \ZZ^y(1)$, et  

\par (ii) $W_{-1}(N_c)$ est parametris\'e 
par le point $c_1$ de $A^x$ et
$( W_{0}/ W_{-1}(N_c))^*$ est param\'etris\'e par le point $c_2$ de 
$(A^*)^{y}$ (ou si on pr\'ef\`ere $ W_{0}/ W_{-2}(N_c)=[\ZZ^x {\buildrel
c_1 \over \longrightarrow } A]$ et $ W_{0}/ W_{-2}(N_c^*)=[\ZZ^y {\buildrel
c_2 \over \longrightarrow } A^*]$). On peut r\'esumer la situation avec 
le diagramme suivant

$$\matrix{
 \ZZ^{x \cdot y}(1)&- & d^* {\PP}& \longrightarrow &
 A^x \times (A^*)^y & \longrightarrow 0 \cr
=&  & \uparrow &  &\uparrow &  \cr
 \ZZ^{x \cdot y} (1)& - & (d^* {\PP})_c&  \longrightarrow & 
[\ZZ^{x} {\buildrel c_1 \over \longrightarrow} A]+
[\ZZ^y {\buildrel c_2 \over \longrightarrow} A^*]  & \longrightarrow 0 \cr}$$

\pn o\`u on a identifie le point $c$ de $A^x \times (A^*)^y(k)$ avec le 1-motif
$[ \ZZ^{x} {\buildrel c_1 \over \longrightarrow} A]+
[\ZZ^y {\buildrel c_2 \over \longrightarrow} A^*]$ et o\`u les traits horizontaux d\'esignent le fait que $ d^* {\PP}$ 
est un $W(1)$-torseur.
\vskip 0.5 true cm

\pn {\bf 3. Etude de ${\rm Gr}_{*}^{W}(\L  {\G}(M))$. }
\vskip 0.5 true cm

\pn 3.1. Soient $M=[X {\buildrel u \over \longrightarrow} G]$ 
un 1-motif de ${\cal M} (k)$ de  gradu\'e 
$\M=X + A + Y(1)$  (i.e. $\M \cong {\rm Gr}_{*}^{W}(M)$), et $\langle M \rangle^\otimes$ la sous-cat\'egorie tannakienne de ${\cal M} (k)$ engendr\'ee par $M$.
Dans cette section
on va calculer de fa\c con explicite les crans ${\rm Gr}_{-1}^{W}$ et 
${\rm Gr}_{-2}^{W}$ de l'alg\`ebre de Lie du groupe de Galois motivique de $M$.

\par On pose 

$$E=W_{-1}( {\underline {\End}}({\M})).$$

\pn  Le dual de ${\M}$ est le 1-motif ${\M}^{\du}=Y^{\du}(-1)+ A^{\du}+X^{\du},$ o\`u $Y^{\du}(-1),A^{\du}$ et $X^{\du}$ sont des motifs purs de poids respectifs 2, 1 et 0.
Puisque ${\underline {\End}}({\M})
\cong {\M}^{\du} \otimes {\M},$ on observe que $E$ est le 1-motif 
scind\'e, de poids $\leq -1,$
de seules composantes non nulles 

$$\eqalign{
E_{-1}&= X^{\du} \otimes A + A^{\du} \otimes Y(1)\cr
E_{-2}&= X^{\du} \otimes Y(1)\cr}$$

\pn de poids respectifs -1 et -2.
D'apr\`es (1.2.1) le 1-motif $A^{\du} \otimes Y(1)$ est isomorphe au 1-motif
$A^* \otimes Y$. De plus, $ A^* \otimes Y(\ok)=(A^{*})^{{\rg} Y}(\ok)$ et 
 $X^{\du} \otimes A(\ok)=A^{{\rg} X^{\du}}(\ok).$

\par La composition des endomorphismes d\'efinit un produit 
$$P:E \otimes E \longrightarrow E$$

\pn sur $E$ de seule composante non nulle le morphisme $E_{-1} \otimes E_{-1} \longrightarrow E_{-2}$ de $\langle M \rangle^\otimes$ qui s'explicite ainsi:

$$E_{-1} \otimes E_{-1} \longrightarrow (X^{\du} \otimes A) \otimes (A^* \otimes Y) 
\longrightarrow {\ZZ}(1) \otimes X^{\du} \otimes Y = E_{-2}$$

\pn o\`u la premi\`ere fl\`eche est la projection de $ E_{-1} \otimes E_{-1}$ sur le
facteur $(X^{\du} \otimes A) \otimes (A^* \otimes Y) $ et la deuxi\`eme fl\`eche
provient du morphisme $P_{{\cal P}}:
A \otimes A^* \longrightarrow \ZZ(1)$ d\'efini par la biextension
 de Poincar\'e ${\cal P}$ de $(A,A^*)$ par $\ZZ(1)$ (cf. (1.3.1)). Ce produit 
$ P : E_{-1} \otimes E_{-1} \longrightarrow E_{-2}$  munit $E$ d'une structure 
d'anneau.
 L'action de $E=W_{-1}( {\underline {\End}}({\M}))$ sur $\M$ est d\'ecrite par le morphisme 
$ E \otimes {\M} \longrightarrow {\M}$ de $\langle M \rangle^\otimes$
 d\'efini par 

 $$\eqalign{
\alpha_1:& (X^{\du} \otimes A) \otimes X  \longrightarrow A \cr
\alpha_2:& (A^* \otimes Y) \otimes A  \longrightarrow Y(1) \cr 
\gamma :& (X^{\du} \otimes Y(1)) \otimes X  \longrightarrow Y(1) \cr}\leqno(3.1.1)$$

\pn o\`u la premi\`ere et la derni\`ere fl\`eche se d\'eduisent
 du morphisme \'evaluation
$ev_{X^{\du}}: X^{\du} \otimes X \longrightarrow \ZZ(0)$, 
tandis que la deuxi\`eme
provient de l'accouplement de Weil $P_{{\cal P}}:
A \otimes A^* \longrightarrow \ZZ(1).$ 
Les morphismes $\alpha_1$ 
et $\gamma$ sont donc des projections.
\pn Le morphisme $(\alpha_1, \alpha_2, \gamma ):  E \otimes {\M} \longrightarrow {\M}$ munit
le 1-motif $\M$ d'une structure de $E$-module.
\vskip 0.5 true cm

\pn 3.2. D'apr\`es (1.3.1) le produit 
$P:E_{-1} \otimes E_{-1} \longrightarrow E_{-2}$
d\'efinit une biextension $\cal B$ de $(E_{-1},E_{-1})$ par 
$E_{-2} .$ Si $d:  E_{-1} \longrightarrow E_{-1} \times E_{-1}$ d\'esigne le
morphisme diagonal de $ E_{-1}$, d'apr\`es le corollaire 2.7
le  $\Sigma - X^{\du} \otimes Y  (1)$-torseur 

$$ d^* {\cal B}=((d^* {\cal B}, -1), ({\cal B} 
\wedge s^*{\cal B}, \xi_{{\cal B} \wedge s^*{\cal B}}), \gamma_{\cal B}, 1 \wedge \nu_{\cal B})$$

\pn  munit le 1-motif $E$ du
 crochet de Lie $[\, ,\,]: E \otimes E \longrightarrow E$, de seule composante non nulle 

$$[\, ,\,]: E_{-1} \otimes E_{-1} \longrightarrow E_{-2}.$$ 

\pn Ce crochet $[\, ,\,]$, qui correspond via (1.3.1) \`a la biextension 
sym\'etris\'ee  ${\cal B} \wedge s^*{\cal B}$,
  est l'antisym\'etris\'e du produit 
$P:E_{-1} \otimes E_{-1} \longrightarrow E_{-2}$ (cf. remarque 1.4).
 Le 1-motif $E$, muni de ce crochet de Lie, 
est une alg\`ebre de Lie $(E, [\, ,\,])$.
\par Puisque le diagramme

%%%%%%%%%%%%%%%%%%%%%%%%%%%%%%%%%%%%%%%%%%%%%%%%%%%%%%%%%%%%%%%%%%%%%%%%%%%%%%
\diagram[small,textflow]
E_{-1} \otimes E_{-1} \otimes X & \rTo^{id_{E_{-1}} \otimes  \alpha_1} &E_{-1} \otimes A & \rTo^{\alpha_2} & Y(1)  \\
 &\rdTo_{[\,,\,] \otimes id_{X}} &  & \ruTo_{\gamma} &  \\
 & & E_{-2} \otimes X & &  \\
\enddiagram 
%%%%%%%%%%%%%%%%%%%%%%%%%%%%%%%%%%%%%%%%%%%%%%%%%%%%%%%%%%%%%%%%%

\pn est commutatif, le morphisme 
$(\alpha_1, \alpha_2, \gamma ):  E \otimes {\M} \longrightarrow {\M}$ est compatible au crochet de Lie
$[\, ,\,]: E_{-1} \otimes E_{-1} \longrightarrow E_{-2}$ et donc
${\M}$ est un $(E, [\, ,\,])$-module de Lie. On a d\'emontr\'e
\vskip 0.5 true cm

\pn 3.3. {\bf Lemme}
\pn {\it 
\par (1) $E=W_{-1}( {\underline {\End}}({\M}))$ est muni d'une 
structure d'alg\`ebre de Lie, $[\, ,\,]: E \otimes E \longrightarrow E$,
 qui \'equivaut via 2.7 \`a la donn\'ee du $\Sigma - X^{\du} \otimes Y  (1)$-torseur 
$ d^* {\cal B}= ((d^* {\cal B}, -1), ({\cal B} \wedge s^*{\cal B},
 \xi_{{\cal B} \wedge s^*{\cal B}}), \gamma_{\cal B}, 1 \wedge \nu_{\cal B})$;
\par (2) ${\M}$ est muni d'une structure de $(E, [\, ,\,])$-module de Lie.}
\vskip 0.5 true cm

\pn 3.4. Remarques: 
\par (1) La donn\'ee des fl\`eches $ \alpha_1 $ et $\gamma$ est \'equivalente 
\`a la donn\'ee des homomorphismes $\gal$-\'equivariants
$\alpha_1: A^{{\rg} X^{\du}}(\ok) \otimes X(\ok)  \longrightarrow A(\ok)$ et 
$\gamma : Y(1)^{{\rg} X^{\du}}(\ok) \otimes X(\ok)  \longrightarrow Y(1)(\ok)$ respectivement.

\par (2) La donn\'ee du morphisme
$\alpha_2:  (A^* \otimes Y) \otimes A \longrightarrow Y(1)$ 
est \'equivalente 
\`a la donn\'ee de l'homomorphismes $\gal$-\'equivariant
$\alpha_2: (A^*)^{{\rg} Y}(\ok) \otimes A(\ok) \longrightarrow Y(1)(\ok)$ 
qui
s'identifie 
\`a ${\rg}(Y)$-copies de l'accouplement de Weil $P_{{\cal P}}:A \otimes A^* \longrightarrow \ZZ(1).$ Par cons\'equent, sur $\ok$ la biextension
 $ \cal B$  est la biextension $\PP$ construite de fa\c con explicite 
dans l'exemple 2.8 et le crochet de Lie $[\, ,\,]: E \otimes E \longrightarrow E$ est d\'ecrit
par les formules (2.8.4) (il suffit de prendre  
 $x={\rg}X^{\du}$ et
$y={\rg}Y$). 

\par (3) Le dual de Cartier de ${\M}$ est le 1-motif scind\'e
${\M}^*= Y^{\du} + A^* + X^{\du}(1).$
 D'apr\`es (1.2.1), on observe que 

$$\eqalign{
{\underline {\End}}({\M}^*) \cong ({\M}^*)^{\du} \otimes {\M}^* & =  
({\M}^{\du}(1))^{\du} \otimes {\M}^{\du}(1)\cr
&= {\M}(-1) \otimes {\M}^{\du}(1) \cr
&= {\M} \otimes {\M}^{\du} \cong {\underline {\End}}({\M}). \cr}$$

\pn Par cons\'equence, le 1-motif $E$ s'identifie avec
$ W_{-1}( {\underline {\End}}({\M}^*))$ et son action 
 $E \otimes {\M}^* \longrightarrow {\M}^*$ sur $ {\M}^*$ est d\'ecrite par les fl\`eches de $\langle M \rangle^\otimes$ 

 $$\eqalign{
\alpha_2^*:& (A^* \otimes Y) \otimes Y^{\du} \longrightarrow A^* \cr
\alpha_1^*:& (X^{\du} \otimes A) \otimes A^*  \longrightarrow X^{\du}(1) \cr 
\gamma^* :& (X^{\du} \otimes Y(1)) \otimes Y^{\du} \longrightarrow X^{\du}(1) \cr}\leqno(3.4.1)$$

\pn o\`u $\alpha_2^*$ et $\gamma^*$ sont des projections tandis que
 $\alpha_1^*$ est donn\'ee par la biextension de Poincar\'e $\cal P$ de $(A,A^*)$.
\vskip 0.5 true cm

\pn 3.5. Soit  $(X,Y^{\du}, A,A^*, v:X \longrightarrow A, 
v^*:Y^{\du} \longrightarrow A^*, \psi:X \otimes Y^{\du}\longrightarrow
(v \times v^*)^*{\cal P})$ le 7-uplet qui d\'efinit le 1-motif 
$M=[X {\buildrel u \over \longrightarrow} G]$ (cf. 1.1). 
Le $\Sigma-X^{\du}\otimes Y(1)$-torseur $d^*{\cal B}$ d\'efini en 3.2, 
permet d'interpreter 
les donn\'ees $v, v^*$ et $\psi$ en termes de points:  
\par Gr\^ace aux morphismes 
$\delta_{ X^{\du}}: \ZZ(0) \longrightarrow X \otimes X^{\du}$ et
$ev_{X}: X \otimes X^{\du} \longrightarrow \ZZ(0)$,
la donn\'ee du morphisme 
$v: X \longrightarrow A$ est \'equivalente \`a la donn\'ee du morphisme 

$$V: \ZZ(0)\,\, {\buildrel \delta_{ X^{\du}} \over \longrightarrow} \,\,
 X \otimes  X^{\du}\, \,
{\buildrel {v \otimes Id_{ X^{\du}}} \over \longrightarrow}\,\,
 A \otimes X^{\du}.$$

\pn De fa\c con analogue, le morphisme $v^*: Y^{\du} \longrightarrow A^*$
 \'equivaut au morphisme

$$V^*: \ZZ(0)\,\, {\buildrel \delta_{Y} \over \longrightarrow} \,\,
Y^{\du} \otimes Y
\,\,{\buildrel {v^* \otimes Id_{ Y}} \over \longrightarrow}\,\,
 A^* \otimes Y.$$

\pn Par cons\'equent se donner les homomorphismes $\gal$-\'equivariants
 $v: X(\ok) \longrightarrow A(\ok)$ et 
$v^*: Y^{\du}(\ok) \longrightarrow A^*(\ok)$ est \'equivalent \`a se donner un point $b=(b_1,b_2)$ de
$E_{-1}(k)= A \otimes X^{\du}(k)+A^* \otimes Y(k).$

\par La trivialisation $\psi$ de l'image r\'eciproque  
 par $(v \times v^*)$ de la biextension de Poincar\'e ${\cal P}$ de $(A,A^*)$
correspond \`a un point $\widetilde b$
dans la fibre de $d^*{\cal B}$ au dessus de $b=(b_1,b_2):$
Le groupe des caract\`eres du tore $X^{\du} \otimes Y(1)$ est le $\gal$-module $X \otimes  Y^{\du}$. 
Fixons un  \'el\'ement $(x,y^{\du})$ de $X \otimes Y^{\du}(\ok)$. 
Par d\'efinition du point $b$,
il existe un \'el\'ement $s$ de $X(\ok)$ et un \'el\'ement
$t $ de $Y^{\du}(\ok)$
tels que $v(x)= \alpha_1(b_1,s) \in A(k)$ et 
$v^*(y^{\du}) =\alpha_2^*(b_2,t)\in A^*(k).$
Notons $i^*_{x,y^{\du}} d^*{\cal B}$ 
la restriction de $d^* {\cal B}$
 par l'inclusion $i_{x,y^{\du}}: \{ (v(x),v^*(y^{\du}) )\}  
 \longrightarrow E_{-1}$ dans $E_{-1}$ de
la sous-vari\'et\'e ab\'elienne engendr\'ee par le point
$ (v(x),v^*(y^{\du}) )$ de $E_{-1}(k)$.
L'image directe 
$(x,y^{\du})_*i^*_{x,y^{\du}} d^*{\cal B}$ 
de $i^*_{x,y^{\du}} d^*{\cal B}$ par le caract\`ere
$(x,y^{\du}):X^{\du} \otimes Y(1)\longrightarrow \ZZ(1)$ est un $\Sigma-\ZZ(1)$-torseur 
sur $ \{ (v(x),v^*(y^{\du}))\}:  $

$$\matrix{ 
(x,y^{\du})_*i^*_{x,y^{\du}} d^*{\cal B}&  \longleftarrow & i^*_{x,y^{\du}} d^*{\cal B} &  \longrightarrow & d^* {\cal B}  \cr
\downarrow & &\downarrow & & \downarrow  \cr 
 \{ (v(x),v^*(y^{\du})) \}    & =  &
 \{ (v(x),v^*(y^{\du})) \}    &  
{\buildrel i_{x,y^{\du}} \over \longrightarrow}   & E_{-1} 
  \cr}$$

\pn D'apr\`es la remarque 3.4 (2), le point $\psi(x,y^{\du})$
correspond \`a un point $(\widetilde b)_{x,y^{\du}}$ de 
\pn $(x,y^{\du})_*i^*_{x,y^{\du}} d^*{\cal B}$ au dessus 
de $ (v(x),v^*(y^{\du}))$. Puisque la connaissance de la famille 
de points $\{(\widetilde b)_{x,y^{\du}}\}_{x,y^{\du}}$ \'equivaut \`a la connaissance d'un point $k$-rationnel $\widetilde b$ de $d^* {\cal B}$ au dessus du point $b$ de $E_{-1}(k), $
 la donn\'ee de la trivialisation
$\psi$ \'equivaut \`a la donn\'ee du point $\widetilde b$ de $d^* {\cal B}$.
\par\noindent On peut donc conclure:
\vskip 0.5 true cm

\pn 3.6. {\bf Lemme}
\pn {\it La donn\'ee du 1-motif $M=[X {\buildrel u \over \longrightarrow} G]$ est \'equivalente aux donn\'ees
\par (a) du gradu\'e $\M=X + A + Y (1)$ de $M$;
\par (b) d'un point $b=(b_1,b_2)$ dans $E_{-1}(k)=X^{\du} 
\otimes A (k) + A^* \otimes Y(k)$; 
\par (c) d'un point $k$-rationnel $\widetilde b$ dans la fibre de $d^* {\cal B}$ au dessus du point $b$ de $E_{-1}$.}
\vskip 0.5 true cm

\pn 3.7. Avant de caract\'eriser la cat\'egorie tannakienne engendr\'ee
par le 1-motif  $M=[X {\buildrel u \over \longrightarrow} G],$ 
on a besoin d'introduire quelques notations: 
\vskip 0.3 true cm

\par - Soit $B$ 
la plus petite sous-vari\'et\'e ab\'elienne (\`a isogenie pr\`es)
de $X^{\du}  \otimes A+A^* \otimes Y$ contenant le point 
$b=(b_1,b_2) \in X^{\du}  \otimes A (k)\times
A^* \otimes Y (k) $.
La restriction $i^*d^* {\cal B}$ de $d^* {\cal B}$ par 
l'inclusion $i: B \longrightarrow E_{-1}$ de $B$ 
dans $E_{-1}$, est un $\Sigma-X^{\du} \otimes Y(1)$-torseur
sur $B$.
\vskip 0.3 true cm

\par - Notons 
$Z_1$ le plus petit sous $\gal$-module de 
$X^{\du} \otimes Y$ tel que le tore $Z_1(1)$, qu'il d\'efini, contienne 
l'image du crochet de Lie $[\, ,\,]: B \otimes B \longrightarrow X^{\du} \otimes Y(1)$.
 L'image directe $p_*i^*d^* {\cal B}$ du $\Sigma-X^{\du} \otimes Y(1)$-torseur
 $i^*d^* {\cal B}$ par la projection
 $p:X^{\du} \otimes Y(1) \longrightarrow (X^{\du} \otimes Y/ Z_1)(1)$ est
 alors un $\Sigma-(X^{\du} \otimes Y/ Z_1)(1)$-torseur trivial sur $B$, 
i.e. $p_*i^*d^* {\cal B}= B \times (X^{\du} \otimes Y/ Z_1)(1)$. On note 

$$\pi: p_*i^*d^* {\cal B} \longrightarrow  (X^{\du} \otimes Y/ Z_1)(1)$$

\pn la projection canonique et  
$s: B \longrightarrow  p_*i^*d^* {\cal B}$
la section canonique. On peut r\'esumer la situation avec le diagramme suivant

$$\matrix{ X^{\du} \otimes Y(1) & =& X^{\du} \otimes Y(1) & {\buildrel p \over \longrightarrow} &(X^{\du} \otimes Y/ Z_1)(1) \cr
| & & | & &\uparrow^\pi \cr
d^* {\cal B} &  \longleftarrow &i^*d^* {\cal B}  & \longrightarrow & p_*i^*d^* {\cal B} \cr
\downarrow & & \downarrow & &~~\downarrow \uparrow^s \cr 
 E_{-1} &  {\buildrel i \over \longleftarrow}   & B & = &B \cr
\downarrow & & \downarrow & &\downarrow \cr
0 & & 0 & & 0 \cr}$$

\pn o\`u les traits verticaux d\'esignent le fait que $d^* {\cal B}$ est un  $X^{\du} \otimes Y(1) $-torseur.
\pn D'apr\`es 3.6 le rel\`evement $u: X \longrightarrow G$ du morphisme 
$v:X \longrightarrow A$ est \'equivalent \`a la donn\'ee du point 
 $\widetilde b$ de $d^* {\cal B}$ au dessus du 
 point $b$. Notons encore $\widetilde b$ les points de $i^*d^* {\cal B}$
et de $p_*i^*d^* {\cal B}$ au dessus du point $b$ de $B$.
\vskip 0.3 true cm

\par - Soit $Z$ le plus petit sous $\gal$-module de $X^{\du} \otimes Y$ 
contenant $Z_1$ et tel que 
le sous-tore $(Z/ Z_1)(1)$ de $(X^{\du} \otimes Y/ Z_1)(1)$ contienne 
$\pi ({\widetilde b}) $. Le tore $Z(1)$ est donc le plus petit 
 sous-tore de  $X^{\du} \otimes Y(1)$ contenant 
$Z_1(1)$ et le point $\pi ({\widetilde b}) .$
\vskip 0.3 true cm

\pn A partir du 1-motif $M=[X {\buildrel u \over \longrightarrow} G]$ on a donc
construit la plus petite sous-alg\`ebre 
de Lie $(B,Z(1), [\, , \,])$ de $( E, [\, , \,])$ qui 
contient toutes les donn\'ees pour d\'efinir 
$\widetilde b$. De plus d'apr\`es 3.3, via la restriction des morphismes (3.1.1) \`a $B$ et $Z(1)$,
$\M$ est un $(B,Z(1),[\,,\,])$-module de Lie.
\vskip 0.5 true cm

\pn 3.8. {\bf Th\'eor\`eme}
\pn {\it  Le foncteur ``prendre le gradu\'e'' 
${\Gr}^W_*:\langle M \rangle^\otimes \longrightarrow
\langle {\M} \rangle^\otimes$ induit une \'equivalence de 
$\langle M \rangle^\otimes$ avec 
 la cat\'egorie des objets de $\langle \M \rangle^\otimes$
munis de l'action de l'alg\`ebre de Lie $(B,Z(1), [\, , \,])$ d\'efinie 
via (3.1.1). }
\vskip 0.5 true cm

\pn PREUVE: Il suffit de prouver qu'\`a partir de l'action de
 $(B,Z(1), [\, , \,])$ sur $\M$ on peut construire un quelconque
$\otimes$-g\'en\'erateur de la cat\'egorie $ \langle M \rangle^\otimes$.
\par Le gradu\'e $\M=X + A + Y(1)$ fournit le 
 quadruplet $(X,Y^{\du},A,A^*),$ o\`u $Y^{\du}$ est le groupe
 des caract\`eres du tore $Y (1)$. 
\par Soient $\{e_i\}_i$ une base de $X(\ok)$ 
et $\{f_j^*\}_j$ une base de $Y^{\du}(\ok)$. 
Choisissons un point $P$ de $B \cap X^{\du} \otimes A (k)$ et un point
$Q$ de $B \cap A^* \otimes Y (k)$ tels que la sous-vari\'et\'e ab\'elienne de
$X^{\du} \otimes A +A^* \otimes Y$ qu'ils engendrent, soit isog\`ene \`a $B$. 
Soient ${\overline v}:X(\ok) \longrightarrow A(\ok)$ et 
${\overline v}^*:Y^{\du}(\ok) \longrightarrow 
A^*(\ok)$ les homomorphismes $\gal$-\'equivariants d\'efinis par 

$$\eqalign{{\overline v}(e_i)&=\alpha_1(P,e_i),\cr 
{\overline v}^*(f_j^*)&=\alpha_2^*(Q,f_j^*).\cr}$$

\par Comme dans 3.7, notons 
$Z_1(1)$ le plus petit sous-tore de $Z(1)$ qui contient
$[B,B]$. Choisissons un point
${\vec q}=(q_1,\dots,q_{{\rm rg}\, (Z/Z_1)})$ 
de $(Z/Z_1)(1)(k)$ tel que les points
$q_1, \dots,q_{{\rm rg}\, (Z/Z_1)} $
 soient multiplicativement ind\'ependants. 

\par Soit $ \Gamma: Z(1)(\ok) \otimes X \otimes Y^{\du}(\ok)
\longrightarrow  \ZZ(1)(\ok)$ l'homomorphisme $\gal$-\'e\-qui\-va\-riant
obtenu \`a partir des fl\`eches 
$ \gamma: (X^{\du} \otimes Y(1)) \otimes X 
\longrightarrow  Y(1)$  et 
$ev_Y: Y \otimes Y^{\du} \longrightarrow \ZZ(0)$. Notons 
${\overline \psi}: X \otimes Y^{\du}(\ok) \longrightarrow {\ZZ}(1)(\ok)$ l'homomorphisme $\gal$-\'equivariant d\'efini par 

$${\overline \psi}(e_i,f_j^*)= \Gamma([P,Q],{\vec q},e_i,f_j^* ).$$ 

\pn D'apr\`es la remarque 3.4 (2), le point ${\overline \psi}(e_i,f_j^*)$ 
fournit 
un point dans la fibre, au dessus de $(e_i,f_j^*),$
de l'image r\'eciproque par $({\overline v} \times {\overline v}^*)$
 de la biextension
 de Poincar\'e ${\cal P}$ de $(A,A^*).$  On a donc que ${\overline \psi}$
est une trivialisation $\gal$-\'equivariant de 
$({\overline v} \times {\overline v}^*)^*{\cal P}.$
\par Les morphismes
${\overline v},
 {\overline v}^*$ et ${\overline \psi}$ d\'efinissent un 
 1-motif
$(X,Y^{\du},A,A^*,{\overline v},{\overline v}^*,{\overline \psi})$, qui
par construction, engendre la m\^eme cat\'egorie tannakienne de $M$.
\vskip 0.5 true cm

\pn 3.9. Remarques:
\par (1) Les points $k$-rationnels $P,Q, {\vec q}$ choisis dans cette d\'emonstration 
permettent seulement de
fixer un $\otimes$-g\'en\'erateur de la cat\'egorie tannakienne
$\langle M \rangle^\otimes$
 et ne permettent pas de retrouver le 1-motif $M$.
En effet, consid\'erons l'exemple suivant: 
\pn soit
$M=[{\ZZ} {\buildrel u \over \longrightarrow } {\GG}_m^3]$
le 1-motif de ${\cal M} (k)$
d\'efini par $u(1)=(q_1,q_2,1)$ avec  $q_1, q_2$ deux \'el\'ements de 
${\GG}_m(k)-\mu_\infty $ multiplicativement ind\'ependants (
$\mu_\infty$ est le groupe des racines de l'unit\'e 
dans $\ok$). 
 On a alors que $Z(1)={\GG}_m^2$ et que la cat\'egorie tannakienne engendr\'ee par $M$ 
est la m\^eme que celles engendr\'ees par les 1-motifs suivants:
\par - $[{\ZZ} {\buildrel u \over \longrightarrow } {\GG}_m^2]$ 
avec $u(1)=(p_1,p_2)$, $p_1$ et $p_2$ deux \'el\'ements de ${\GG}_m(k)-\mu_\infty$ 
multiplicativement ind\'ependants qui appartiennent au groupe multiplicatif engendr\'e par $q_1$ et $q_2$;
\par - $[{\ZZ}^4 {\buildrel u \over \longrightarrow } {\GG}_m]$ avec $u(1,0,0,0)=r_1$,
 $u(0,1,0,0)=r_1^3$, $u(0,0,1,0)=1$
\pn  $u(0,0,0,1)=r_2$, $r_1$ et $r_2$ deux \'el\'ements de ${\GG}_m(k)-\mu_\infty$
 multiplicativement in\-d\'e\-pen\-dants qui appartiennent au groupe multiplicatif 
engendr\'e par $q_1$ et $q_2$.

\par (2) Gr\^ace \`a ce th\'eor\`eme, 
on peut conclure que $(B,Z(1), [\, , \,])$ est la plus petite 
sous-alg\`ebre 
de Lie de $( E, [\, , \,])$
qui a la propri\'et\'e de caract\'eriser la cat\'egorie  
$\langle M \rangle^\otimes$ au moyen de son action sur le gradu\'e $\M$ de $M$.
\vskip 0.5 true cm

\pn 3.10. {\bf Corollaire}
\pn {\it Le gradu\'e par le poids ${\rm Gr}_{*}^{W}(\L W_{-1} {\G}(M))$
 de l'alg\`ebre de Lie du radical unipotent
du groupe de Galois motivique de $M$ est l'alg\`ebre de Lie $(B,Z(1), [\, , \,]).$ 
\pn En particulier, ${\rm Gr}_{-1}^{W}(\L {\G}(M))=(B, [\, , \,]) $ et 
$W_{-2}(\L {\G}(M))=(Z(1), [\, , \,]). $ }
\vskip 0.5 true cm

\pn PREUVE: Le foncteur ``prendre le gradu\'e'' 
${\Gr}^W_*:\langle M \rangle^\otimes \longrightarrow
\langle {\M} \rangle^\otimes$ fournit l'inclu\-sion de sch\'emas motiviques 
 
$$ {\G}(\M) \longrightarrow {\Gr}^W_* ({\G}(M))$$

\pn qui traduit le fait que le groupe de Galois motivique
 du gradu\'e $\M$ est le quotient 
${\rm Gr}_{0}^{W}$ du groupe de Galois motivique de $M$.
\pn Si on applique le th\'eor\`eme 8.17 [D90] \`a ce foncteur ${\Gr}^W_*$, 
on obtient que la cat\'egorie  $ \langle M \rangle^\otimes$ est \'equivalente \`a la
cat\'egorie
des objets de  $\langle {\M} \rangle^\otimes$ munis d'une action de  
${\rm Gr}_{*}^{W}( {\G}(M))$ factorisant l'action de 
${\G}(\M)={\rm Gr}_{0}^{W}( {\G}(M))$. 
En passant aux alg\`ebres de Lie, la cat\'egorie
$ \langle M \rangle^\otimes$ est donc \'equivalente \`a la cat\'egorie
des objets de  $\langle {\M} \rangle^\otimes$ munis d'une action de  

$${\rm Gr}_{-1}^{W}(\L {\G}(M))+W_{-2}(\L {\G}(M)).$$

\pn D'un autre c\^ot\'e, d'apr\`es le th\'eor\`eme 3.8 la cat\'egorie 
$\langle {M} \rangle^\otimes$ est \'equivalente 
 \`a la cat\'egorie
des objets de  $\langle {\M} \rangle^\otimes$ munis de l'action de  
l'alg\`ebre de Lie 

$$(B,Z(1),[\,,\,])$$

\pn d\'ecrite en (3.1.1). La
 propri\'et\'e universelle (1.5.2) de ${\G}(M)$ nous fournit alors
 l'inclusion

$${\rm Gr}_{-1}^{W}(\L {\G}(M))+W_{-2}(\L {\G}(M))\longrightarrow
 (B,Z(1),[\,,\,]).$$

\pn Mais d'apr\`es la remarque 3.9 (2), $ (B,Z(1),[\,,\,])$ est la plus 
petite sous-alg\`ebre de Lie qui a la propri\'et\'e de caract\'eriser la cat\'egorie 
$\langle {M} \rangle^\otimes$ au moyen de son action sur le gradu\'e
 $\M$, et donc on peut conclure que
 cette inclusion est en r\'ealit\'e un isomorphisme, i.e.

$${\rm Gr}_{-1}^{W}(\L {\G}(M))+W_{-2}(\L {\G}(M)) \cong
 (B,Z(1),[\,,\,]).$$

\pn 3.11. Remarques:
\par (1) Si le 1-motif $M$ est quasi-d\'eficient, i.e. si $W_{-1}(\G (M))$ est ab\'elien 
(cf. [B02]), la vari\'et\'e ab\'elienne $B$
correspond \`a la vari\'et\'e ab\'elienne isotrope maximale introduite par D. Bertrand dans la preuve
du th\'eor\`eme 1 [Be98].
\par (2) D'apr\`es l'analogue motivique de 1.4 et (1.3.1) [B02], le sous-tore $Z_1(1)$ de 
$ W_{-2}(\L {\G}(M))$ est le groupe d\'eriv\'e de 
$ W_{-1}(\L {\G}(M)).$
\vskip 0.5 true cm

\pn {\bf 4. Preuve du th\'eor\`eme principal}
\vskip 0.5 true cm

\pn 4.1. Gr\^ace aux r\'esultats du paragraphe 3, on peut finalement prouver 
que le radical unipotent $W_{-1}(\L {\G}(M))$
 de l'alg\`ebre de Lie du 
groupe de Galois motivique de $M$ est la vari\'et\'e 
semi-ab\'elienne extension de $B$ par $Z(1)$ d\'efinie
par l'action adjointe de l'alg\`ebre de Lie $(B,Z(1), [\, , \,])$ sur $B+Z(1).$
Pour faire ceci, on doit d'abord g\'en\'eraliser la preuve du th\'eor\`eme 3.8.
\par Soit $(C+W(1),[\,,\,])$ une alg\`ebre de Lie de ${\cal M}(k)$, o\`u 
$C$ est une vari\'et\'e ab\'elienne et $W(1)$ un tore. D'apr\`es la proposition 2.5,
le crocher de Lie $[\,,\,]$ fournit un $\Sigma-W(1)$-torseur $\cal D$ sur $C$.
\vskip 0.5 true cm

\pn 4.2. {\bf Lemme}
\pn {\it Soient $c$ un point de $C(k)$ et $\widetilde c$ un point $k$-rationnel 
dans la 
fibre de $\cal D$ au dessus de $c$. 
\pn Une action fid\`ele de l'alg\`ebre de Lie $(C+W(1),[\, , \,])$ 
sur un 1-motif $\widetilde N$ de filtration par le poids scind\'ee, fournit
un 1-motif $N$ de gradu\'e $\widetilde N$.}
\vskip 0.5 true cm

\pn PREUVE: Soit ${\widetilde N}=X+A+Y(1)$ un 1-motif de ${\cal M}(k)$
de filtration par le poids scind\'ee.
\pn Puisque l'action est fid\`ele, le 1-motif
$ C+W(1)$ s'injecte dans le motif $ {\underline {\End}}({\widetilde N}) \cong {\widetilde N}^{\du} 
\otimes {\widetilde N}$
et donc $C$ est une sous-vari\'et\'e ab\'elienne de $ X^{\du} \otimes A+A^* \otimes Y$ et 
$ W(1)$ est un sous-tore de $ X^{\du} \otimes Y(1)$. 
\par L'action de 
l'alg\`ebre de Lie $(C+W(1),[\, , \,])$ sur $\widetilde N$ est d\'efinie par la restriction \`a $C$ et \`a $W(1)$ 
des fl\`eches de $\langle {\widetilde N} \rangle^\otimes$

$$\eqalign{
\alpha_1 &: (X^{\du} \otimes A) \otimes X \longrightarrow A \cr
\alpha_2 &: (A^* \otimes  Y)   \otimes A \longrightarrow Y(1)\cr
\gamma &: (X^{\du} \otimes Y(1)) \otimes X  \longrightarrow Y(1).\cr}$$

\pn 
 o\`u la premi\`ere et la derni\`ere fl\`eche se d\'eduisent du morphisme \'evaluation
$ev_{X^{\du}}: X^{\du} \otimes X \longrightarrow \ZZ(0)$, tandis que la deuxi\`eme
provient de l'accouplement de Weil $P_{{\cal P}}:A \otimes A^* \longrightarrow \ZZ(1).$ 
Les morphismes $\alpha_1$ et $\gamma$ sont donc des projections. 
\pn Le dual de Cartier de ${\widetilde N}$ est le 1-motif 
 ${\widetilde N}^*= Y^{\du} + A^* + X^{\du}(1)$. Comme on l'a d\'ej\`a 
remarqu\'e en 3.4 (3), les motifs
 $ {\underline {\End}}({\widetilde N}^*)$ et $ {\underline {\End}}({\widetilde N})$ sont isomorphes, et donc
l'alg\`ebre de Lie 
$(C+W(1),[\, , \,])$ agit aussi sur ${\widetilde N}^*$ via les fl\`eches 

 $$\eqalign{
\alpha_2^*:& (A^* \otimes  Y) \otimes Y^{\du} \longrightarrow A^* \cr
\alpha_1^*:& (X^{\du} \otimes A) \otimes A^*  \longrightarrow X^{\du}(1) \cr 
\gamma^* :& (X^{\du} \otimes Y(1)) \otimes Y^{\du} \longrightarrow X^{\du}(1) \cr}$$

\pn o\`u $\alpha_2^*$ et $\gamma^*$ sont des projections tandis que  
$\alpha_1^*$ est donn\'ee par la biextension de Poincar\'e $\cal P$ de $(A,A^*)$.

\par Soient $c_1$ et $c_2$ les projections du point $c$ 
dans $A \otimes X^{\du}(k)$ et $A^* \otimes Y(k)$ respectivement et
notons 
$V: \ZZ(0)(\ok) \longrightarrow A \otimes X^{\du}(\ok)$ et 
$V^*: \ZZ(0)(\ok) \longrightarrow A^* \otimes Y(\ok)$ les homomorphismes 
$\gal$-\'equivariants d\'efinis par $c_1$ und $c_2$. 
Gr\^ace aux fl\`eches \'evaluations 
$ev_{X}: X \otimes X^{\du} \longrightarrow \ZZ(0)$ et 
$ev_{Y}: Y \otimes Y^{\du} \longrightarrow \ZZ(0)$, \`a partir de $V$ et $V^*$ 
on obtient des 
 homomorphismes $\gal$-\'equivariants $v: X(\ok) \longrightarrow  A(\ok)$ et 
$v^*: Y^{\du}(\ok) \longrightarrow  A^*(\ok)$.
 On a donc construit 
le 6-uplet $(X,Y^{\du}, A,A^*,v,v^*).$
\par Prouvons maintenant que la donn\'ee du point $\widetilde c$
 correspond \`a une trivialisation $\psi$ de 
de l'image r\'eciproque par $(v \times v^*)$ de la biextension
 de Poincar\'e ${\cal P}$ de $(A,A^*)$.
\pn Le point $\widetilde c$ de $ {\cal D}$
est d\'efini par le point $c$ de $C(k)$ et par un point ${\vec z}$ de $W(1)(k)$.
Soit $ \Gamma: W(1)(\ok) \otimes X \otimes Y^{\du}(\ok)
\longrightarrow  \ZZ(1)(\ok)$ l'homomorphisme $\gal$-\'equivariant
obtenu \`a partir des fl\`eches 
$ \gamma: W(1) \otimes X 
\longrightarrow  Y(1)$  et 
$ev_Y: Y \otimes Y^{\du} \longrightarrow \ZZ(0)$. Notons 
$\psi:X \otimes Y^{\du}(\ok) \longrightarrow {\ZZ}(1)(\ok)$
 l'homomorphisme $\gal$-\'equivariant d\'efini par 

$$\psi(s,t)= \Gamma({\vec z},s,t )$$ 

\pn pour chaque \'el\'ement $(s,t)$ de $X \otimes Y^{\du}(\ok)$.
D'apr\`es la remarque 3.4 (2), le point $\psi(s,t)$ fournit 
un point dans la fibre, au dessus de $(s,t),$ de l'image r\'eciproque 
par $(v \times v^*)$ de la biextension de Poincar\'e ${\cal P}$ de $(A,A^*)$,
i.e. $\psi $
est une trivialisation $\gal$-\'equivariant de 
$(v \times v^*)^*{\cal P}.$
On a donc construit un 1-motif
$N=(X,Y^{\du},A,A^*, v,v^*,\psi)$ dont le gradu\'e est $\widetilde N$.
\vskip 0.5 true cm

\pn 4.3. A partir de maintenant on va reprendre les notations de 3.7.
D'apr\`es le corollaire 3.10, ${\Gr}^W_{-1}(\L {\G}(M))$
est la plus petite
 sous-vari\'et\'e ab\'elienne $B$ de 
$X^{\du} \otimes A + A^* \otimes Y $ contenant le point $k$-rationnel 
$b=(b_1,b_2)$ et $W_{-2}(\L {\G}(M))$ est
le plus petit sous-tore
$Z(1)$ de $X^{\du} \otimes Y(1)$ contenant les donn\'ees
nec\'essaires pour d\'efinir le point $k$-rationnel $\tilde b$ dans la fibre de $d^* {\cal B}$
 au dessus du point $b$.
\vskip 0.5 true cm

\pn 4.4. PREUVE du th\'eor\`eme principal: 
\pn D'apr\`es 3.10 il est clair que la suite (0.1.1) est exacte.
Les morphismes (1.5.1)
d\'efinissent une action du groupe ${\G}(M)$ sur chaque \'el\'ement
de $\langle M \rangle^\otimes$. Par passage aux Ind-objets, on obtient une action 
de ${\G}(M)$ sur tout sch\'ema affine motivique en $\langle M \rangle^\otimes$.
En particulier, ${\G}(M)$ agit sur lui-m\^eme 
par automorphismes int\'erieurs (cf. [D89] 6.1) et donc
$(B,Z(1), [\, , \,])$ agit sur 
$B+Z(1)$ par action adjointe.  
Si on applique le lemme 4.2 \`a cette action et au point $b$ de $B(k)$ on obtient le
1-motif $[0 \longrightarrow W_{-1}(\L {\G}(M))].$
\vskip 0.5 true cm

\pn 4.5. Remarque: Le dual de Cartier du 1-motif 
$[0 \longrightarrow W_{-1}(\L {\G}(M))]$ peut se construire de fa\c con explicite: 
\pn Le dual de Cartier de $B+Z(1)$ est le 1-motif $Z^{\du}+B^*$, o\`u
$Z^{\du}$ est le sous $\gal$-module de $X \otimes Y^{\du}$ qui 
correspond au tore $Z(1)$ et $B^*$ est la sous-vari\'et\'e ab\'elienne de $X \otimes A^* + A \otimes Y^{\du}$ duale de $B$.
Puisque ${\underline {\End}}(B+Z(1)) \cong {\underline {\End}}(Z^{\du}+B^* )$, 
l'alg\`ebre de Lie $(B,Z(1), [\, , \,])$ agit aussi sur $Z^{\du}+B^*$. Cette
 action duale a pour seule composante non nulle la restriction
 $ \beta^*_{\vert B, Z^{\du}} : B \otimes Z^{\du} \longrightarrow B^*$
\`a $B$, $B^*$ et $Z^{\du}$ de la fl\`eche

$$ \beta^*: (X^{\du} \otimes A + A^* \otimes Y)
 \otimes X \otimes Y^{\du} \longrightarrow A \otimes Y^{\du}+ X \otimes A^* $$

\pn qui est la compos\'ee des fl\`eches 

$$\eqalign{\beta^*_1 & :(X^{\du} \otimes A) \otimes X \otimes Y^{\du} 
\longrightarrow A \otimes Y^{\du} \cr
\beta^*_2 & : (A^* \otimes Y) \otimes X \otimes Y^{\du} 
\longrightarrow X \otimes A^*}$$

\pn o\`u $\beta^*_1$ s'identifie \`a ${\rg}(Y^{\du})$-copies de la projection 
$\alpha_1:(X^{\du} \otimes A) \otimes X \longrightarrow A$ introduite en (3.1.1) et
$\beta^*_2$ s'identifie \`a  ${\rg}(X)$-copies de la projection 
$\alpha_2^*:(A^* \otimes Y) \otimes Y^{\du} \longrightarrow  A^*$
 introduite en (3.4.1). On a donc que aussi

$$ \beta^*=(\beta^*_1,\beta^*_2):(A^{{\rg} X^{\du}} \otimes X)^{{\rm rg}\,Y^{\du}} +
 ((A^{*})^{{\rg} Y} \otimes Y^{\du})^{{\rm rg}\,X} \longrightarrow A^{{\rm rg}\,Y^{\du}} +
 (A^*)^{{\rm rg}\,X}$$

\pn  est une projection:
c'est la projection qui repr\'esente l'action adjointe de 
l'alg\`ebre de Lie $(B,Z(1), [\, , \,])$ sur $Z^{\du}+B^*$. A partir de cette 
fl\`eche $\beta^*_{\vert B, Z^{\du}}$, on d\'efinit l'ho\-mo\-mor\-phisme $\gal$-\'equivariant

$$\eqalign{{\cal V}: Z^{\du}(\ok) &\longrightarrow B^*(\ok)\cr
 z & \longmapsto \beta^*_{\vert B, Z^{\du}}(b,z)=(\beta^*_1(b_1,z),\beta^*_2(b_2,z))\cr} $$ 

\pn Par construction, le 1-motif $[{\cal V}:Z^{\du} \longrightarrow B^*]$ 
est le dual de Cartier du 1-motif $[0 \longrightarrow W_{-1}(\L {\G}(M))]$.
\vskip 0.5 true cm

\pn {\bf Bibliographie}
\vskip 0.5 true cm

\par\noindent
[B02] C. Bertolin, {\it The Mumford-Tate group of 1-motives,} Ann. de l'Inst. Fourier 52 (2002).

\par\noindent
[Be98] D. Bertrand , {\it Relative splitting of one-motives},
Number Theory (Ti\-ru\-chi\-ra\-palli, 1996), Contemp. Math. 210 Amer. Math.
Soc. (1998).

\par\noindent
[Br83] L. Breen, {\it Fonctions th\^eta et th\'eor\`eme du cube}, LN 980 (1983).

\par\noindent
[By83] J.-L. Brylinski, {\it 1-motifs et formes automorphes
 (th\'eorie arithm\'etique des domaines 
des Siegel)}, Pub. Math. Univ. Paris VII 15 (1983).

\par\noindent
[D75] P. Deligne, {\it Th\' eorie de Hodge III}, Pub. Math. de
 l'I.H.E.S. 44 (1975).

\par\noindent
[D89] P. Deligne, {\it Le groupe fondamental de la droite 
projective moins trois points}, Galois group over $\QQ$, Math. Sci. Res. Inst. Pub. 16 (1989).

\par\noindent
[D90] P. Deligne, {\it Cat\'egories tannakiennes}, The Grothendieck Festschrift II,
Birk\-h\"au\-ser 87 (1990).

\par\noindent
[D01] P. Deligne, lettre \`a l'auteur (2001).

\par\noindent
[J90] U. Jannsen, {\it Mixed motives and algebraic K-theory}, LN 1400 (1990).

\par\noindent
[M94] J.S. Milne, {\it Motives over finite fields}, Motives, Proc. of Symp. in Pure Math. 55 (1994).

\par\noindent
[R94] M. Raynaud, {\it 1-motifs et monodromie g\'eom\'etrique}, Ast\'erisque 223 (1994).

\par\noindent
[S72] N. Saavedra Rivano, {\it Cat\' egories tannakiennes}, LN 265 (1972).

\pn SGA7 I: Groupes de Monodromie en G\'eom\'etrie Alg\'ebrique, LN 288 (1972).
\vskip 2 true cm

\pn \hskip 1 true cm Bertolin Cristiana
\pn \hskip 1 true cm NWF I - Mathematik
\pn \hskip 1 true cm Universit\"at Regensburg
\pn \hskip 1 true cm D-93040 Regensburg
\pn \hskip 1 true cm Germany
\vskip 0.2 true cm
\pn \hskip 1 true cm Email: cristiana.bertolin@mathematik.uni-regensburg.de

\end